\documentclass[11pt]{article}
\usepackage[dvips]{graphics,color}
\usepackage {pstcol}

\newpsobject{showgrid}{psgrid}{subgriddiv=1,griddots=10,gridlabels=6pt}

\usepackage[dvips]{graphics}
\newcommand{\nc}{\newcommand}
\renewcommand{\sp}{\vskip1ex}
\nc{\z}{\zeta}
\nc{\ps}{\psi}
\nc{\inv}{^{-1}}
\nc{\ep}{\varepsilon}
\nc{\e}{\epsilon}
\def\ra{\rightarrow}

\renewcommand{\l}{\ell}
\nc{\si}{\sigma}
\renewcommand{\u}{u_{c'}}
\def\iy{\infty}

\def\be{\begin{equation}}
\def\ee{\end{equation}}
\def\ba{\begin{eqnarray*}}
\def\ea{\end{eqnarray*}}
\def\bae{\begin{eqnarray}}
\def\eae{\end{eqnarray}}
\def\bc{\begin{center}}
\def\ec{\end{center}}
\def\ov{\over}
\def\al{\alpha}
\def\s{\sigma}
\def \tS{\tilde S}

\def\vp{\varphi}

\nc{\noi}{\noindent}

\def\la{\lambda}
\def\pr{\textrm{Prob}}

\def\cd{\cdots}

\def\h{\Delta h}
\def\dl{\delta}

\def\t{\tau}
\def\F{\mathcal{F}}
\nc{\ga}{\gamma}
\def\ph{\varphi}

\def\convd{\buildrel d \over \longrightarrow}

\begin{document}
\title{
Limit Theorems for Height Fluctuations\\
in a\\  Class of
Discrete Space and Time Growth Models \\
}
\date{October 25, 2000}
\author{
Janko Gravner\\
Department of Mathematics\\
University of California\\
Davis, CA 95616, USA\\
e-mail: gravner@math.ucdavis.edu\\
 \\
Craig A.~Tracy\\
Department of Mathematics\\
Institute of Theoretical Dynamics\\
University of California\\
Davis, CA 95616, USA\\
e-mail: tracy@itd.ucdavis.edu\\
  \\
Harold Widom\\
Department of Mathematics\\
University of California\\
Santa Cruz, CA 95064, USA\\
e-mail: widom@math.ucsc.edu
}
\maketitle
\vspace{-8ex}
\begin{abstract}
We introduce a class of one-dimensional discrete space-discrete time
stochastic growth models described by a  height function $h_t(x)$ with
corner initialization. We prove, with one exception,
that the limiting distribution function of $h_t(x)$ (suitably
centered and normalized)   equals
a Fredholm determinant previously encountered in
random matrix theory.  In particular, in
the universal regime of large $x$ and large $t$
the limiting distribution is the Fredholm determinant
with Airy kernel.   In the exceptional case,  called the
critical regime, the limiting distribution seems not to have
previously occurred.
The proofs use
the dual RSK algorithm, Gessel's theorem, the Borodin-Okounkov
identity and a novel,  rigorous saddle point
analysis.
In the fixed $x$, large $t$ regime,
we find
a Brownian motion representation. This model is equilvalent
to the Sepp{\"a}l{\"a}inen-Johansson model.  Hence some
of our results are not new, but the proofs are.
\end{abstract}


\tableofcontents

\section{Introduction}
\setcounter{equation}{0}
Growth processes have been extensively studied by mathematicians
and physicists for many years (see, e.g., \cite{gravnerGriffeath,
lassig, sepp} and references therein), but it was only recently that
K.~Johansson~\cite{johanssonShape} proved that
  the \textit{fluctuations}
of the limiting shape in
a class of growth models are described by certain distribution functions
first appearing in random matrix theory (RMT)~\cite{tw1,tw2}.
Further work
by Johansson~\cite{johanssonDiscrete}, Pr\"ahofer and
Spohn~\cite{praehofer1,
praehofer2} and Baik and Rains~\cite{baikRainsPNG} strongly suggests the
universal nature of these RMT distribution functions.
 These developments are part
of the recent activity relating Robinson-Schensted-Knuth (RSK) type
problems
of combinatorial probability to
the distribution functions of RMT, see
e.g.~\cite{adler1, aldous2, bdj, baikRainsAsy,baikRainsSym,
borodin3,  itw1, itw2, kuperberg, vanM, okounkov, tw4, tw5}.

\begin{figure}
\bc
\resizebox{10cm}{10cm}{\includegraphics{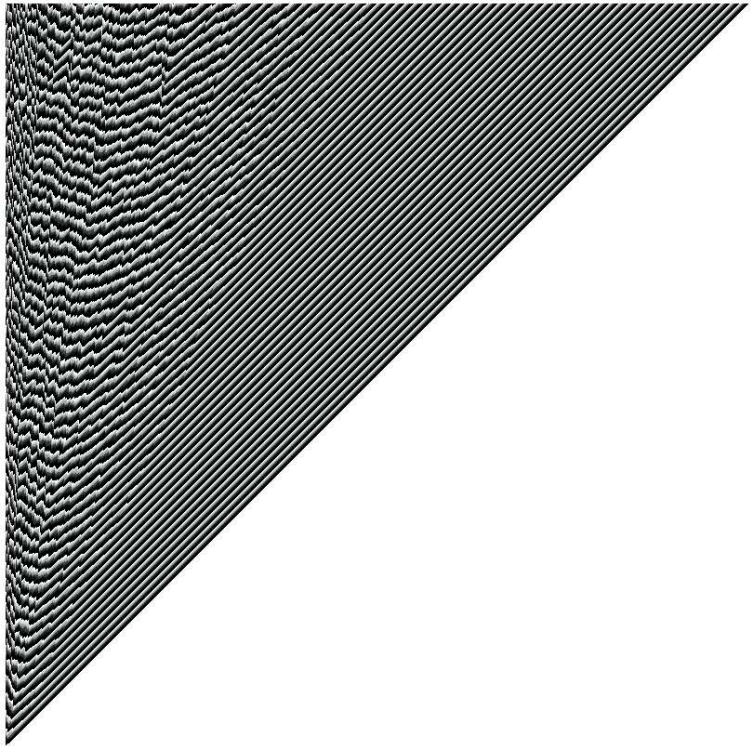}}
\caption{\label{fluctuationsFigs}  ODB Simulations: The contour lines 
of $h_t(x)$ are drawn in  space-time $(x,t)$. In these simulations $0\le t,x\le
800$
and $p=1/2$.}
\ec
\end{figure}

In this paper we analyze a class of one-dimensional discrete
space-discrete
time stochastic growth models,
called \textit{oriented digital boiling}~\cite{gravner1, gravner2}.
  Digital boiling dynamics is
a cellular automaton that models an excitable medium in
the presence of persistent random spontaneous excitation.
Alternatively, digital
boiling models represent contour (constant height) lines for one
of the simplest models for growing a connected interface.
It is this latter point of view we adopt here; that is, we introduce
a height function $h_t(x)$ that characterizes the state of the system.
Fig.~\ref{fluctuationsFigs}  illustrates the height fluctuations
in oriented digital boiling.

We shall derive various limit theorems for $h_t(x)$.
We find four limiting regimes:
\begin{enumerate}
\item[1.] \textit{GUE Universal Regime}: $x\ra\iy$, $t\ra\iy$ such that
$p_c:=1-x/t$ is fixed and $p<p_c$.
\item[2.] \textit{Critical Regime}: $x\ra\iy$, $t\ra\iy$ such that
 $p_c:=1-x/t$ is fixed and $p\sim p_c$.
 \item[3.] \textit{Deterministic Regime}: $x\ra\iy$, $t\ra\iy$ such that
 $p_c:=1-x/t$ is fixed and $p>p_c$.
\item[4.] \textit{Finite $x$ GUE Regime}: Fixed $x$ and $t\ra\iy$.
\end{enumerate}

The limit theorems are stated at the beginning of \S3.  Here is an
outline of how
they are obtained.
First we show that $h_t(x)$ satisfies a last
passage property, i.e.\ it equals the maximum over a certain class of
paths
in space-time. Then applying the
dual RSK algorithm~\cite{knuth, stanley},
we obtain a reformulation of the problem in terms of Young tableaux.
This is
followed by an application of a theorem of Gessel~\cite{gessel} (see
also \cite{tw5})
which gives a Toeplitz determinant representation for the distribution
function
for $h_t(x)$. An identity of Borodin and Okounkov~\cite{borodin2}
expresses the Toeplitz
determinant in terms of the Fredholm determinant of an infinite matrix.
Finally we
use a saddle point analysis (steepest descent) to determinine the
limiting behavior of
the entries, and therefore the Fredholm determinant, of the infinite
matrix.

Along the way
we identify\footnote{A referee points out that this identification can be
made at the very beginning.} ODB with a first-passage
percolation model of
Sepp{\"a}l{\"a}inen~\cite{sepp} whose limit law in the universal
regime  was determined by Johansson~\cite{johanssonDiscrete}.
Thus we could have used  the analysis
in~\cite{johanssonDiscrete}
to establish our limit law in the universal regime, or alternatively
used Riemann-Hilbert
methods \cite{bdj, diz, dz}, to investigate the Toeplitz determinant
asymptotics.
But the method we present is in our opinion more straightforward and
technically simpler
than these, and it is very general. (The Fredholm determinant is easier
to handle than the Toeplitz determinant, even though they are
essentially equal.)
Also, our analysis permits a nice conceptual
understanding of the various limiting regimes.  For example, the
universal regime is
characterized by the coalescence of two saddle points; and the emergence
of
the \textrm{Airy kernel} is related to the well-known appearance of Airy
functions
in such a saddle point analysis~\cite{chester}.

Even in this simpler approach there are technical details to work out
after the
saddle point analysis gives us the answer. For example in the
universal regime we need uniform estimates on the entries of the
infinite matrix
in order to show that the matrix scales in trace norm to the Airy
kernel. These details
 are given completely only for
this regime.

In Regime 4 we give an independent
proof
that the suitably centered and normalized $h_t(x)$
 has a limiting distribution.  The proof proceeds
through the introduction of a certain Brownian motion functional.
This  leads to some apparently new
identities for $n$-dimensional Brownian motion; see (\ref{distrFnGUE})
below.

The initial conditions are corner initialization.  Due to the fact
there is no
known symmetry theorem for the \textit{dual} RSK algorithm~\cite{knuth,
stanley},
we are unable to prove limit theorems with different initial conditions,
e.g.\
growth from a flat substrate.  From work of Baik and
Rains~\cite{baikRainsSym}
and Pr\"ahofer and Spohn~\cite{praehofer1, praehofer2}, it is natural to
conjecture
that the limiting distribution is now of  GOE symmetry and hence given
by the analogous distribution function in the GOE case~\cite{tw2}.

The table of contents provides a detailed description of the
organization of
this paper.

\section{Growth Models and Increasing Paths}
\setcounter{equation}{0}
In this section we introduce three classes of
discrete space and discrete time stochastic growth models.
Each of these models will have an equivalent path
description, but only
for one of these models are we able
to prove limit theorems.
Nevertheless, we believe it is useful to place this
``solvable'' case in a larger context.

We assume that the occupied set of our growth models can be
described by a \textit{height function}
$h_t:\textbf{Z}_+\ra\textbf{Z}_+\cup
\{-\iy,\iy\}$, where
$\textbf{Z}_+$
is the set of nonnegative
integers.   Here, time $t=0,1,2,\ldots$ proceeds in discrete steps.
The \textit{occupied set} at time $t$ is thus given by
\[ \eta_t=\left\{(x,y)\in\textbf{Z}_+\times\textbf{Z}_+: y\le h_t(x)
\right\}. \]
In the models below we use the following one-dimensional
neighborhood: $(x+\mathcal{N})=\{x-1,x\}$ (\textit{the oriented case})
and assume \textit{corner initialization},
\be h_0(x)=\left\{\begin{array}{rl}
                0, &\textrm{if}\> x=0,\\
                -\iy, &\textrm{otherwise}.
        \end{array}\right. \label{initial}\ee

\subsection{Oriented Digital Boiling}
The first class of growth rules  we call
\textit{oriented digital boiling} (ODB)~\cite{gravner1,
gravner2}.\footnote{For
spatial dimensions greater than one, visual features of this dynamics
resemble bubble formation, growth and annihilation in a boiling
liquid, hence the process is called \textit{digital boiling} and
 \textit{oriented} refers to the choice of neighborhood $\mathcal{N}$,
see Fig.~2 in \cite{gravner1} or Feb.~12, 1996 Recipe of
\cite{griffeath2}.}
The rules for ODB are
\begin{enumerate}
\item[(1)] $h_t(x)\le h_{t+1}(x)$ for all $x$ and $t$.
\item[(2)] If $h_t(x-1)>h_t(x)$, then
$h_{t+1}(x)=h_t(x-1)$.
\item[(3)] Otherwise,  then independently
of the other sites and other times, $h_{t+1}(x)=h_t(x)+1$
with probability $p$.  (With probability $1-p$,
we have $h_{t+1}(x)=h_t(x)$.)
\end{enumerate}
It follows from these rules that for every $x$ and $t$,
$ h_t(x-1)\le h_t(x)+1$.

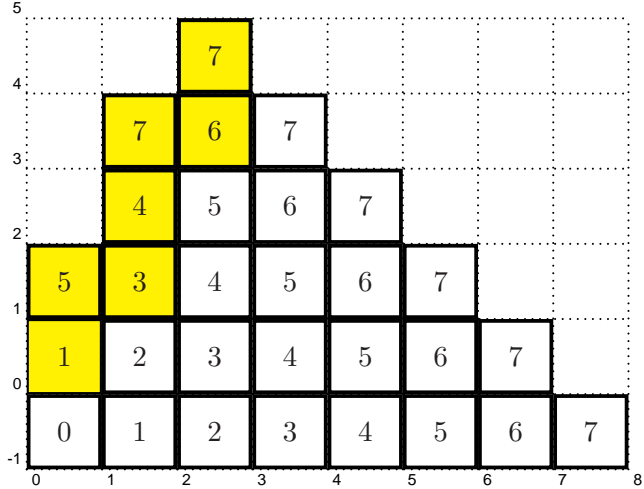
\begin{figure}
\bc
\begin{pspicture}(0,-1)(8,5)\showgrid
\psframe[linestyle=solid,linewidth=.5mm](0,-1)(1,0)
\rput(.5,-.5){0}
\psframe[linestyle=solid,linewidth=.5mm](1,-1)(2,0)
\psframe[linestyle=solid,fillstyle=solid,fillcolor=yellow](0,0)(1,1)
\rput(.5,.5){1}
\rput(1.5,-.5){1}
\psframe[linestyle=solid,linewidth=.5mm](2,-1)(3,0)
\psframe[linestyle=solid,linewidth=.5mm](1,0)(2,1)
\rput(1.5,.5){2}
\rput(2.5,-.5){2}
\psframe[fillstyle=solid,fillcolor=yellow](1,1)(2,2)
\psframe[linestyle=solid,linewidth=.5mm](1,1)(2,2)
\psframe[linestyle=solid,linewidth=.5mm](2,0)(3,1)
\psframe[linestyle=solid,linewidth=.5mm](3,-1)(4,0)
\rput(1.5,1.5){3}
\rput(2.5,.5){3}
\rput(3.5,-.5){3}
\psframe[fillstyle=solid,fillcolor=yellow](1,2)(2,3)
\psframe[linestyle=solid,linewidth=.5mm](1,2)(2,3)
\psframe[linestyle=solid,linewidth=.5mm](2,1)(3,2)
\psframe[linestyle=solid,linewidth=.5mm](3,0)(4,1)
\psframe[linestyle=solid,linewidth=.5mm](4,-1)(5,0)
\rput(1.5,2.5){4}
\rput(2.5,1.5){4}
\rput(3.5,.5){4}
\rput(4.5,-.5){4}
\psframe[fillstyle=solid,fillcolor=yellow](0,1)(1,2)
\psframe[linestyle=solid,linewidth=.5mm](0,1)(1,2)
\psframe[linestyle=solid,linewidth=.5mm](2,2)(3,3)
\psframe[linestyle=solid,linewidth=.5mm](3,1)(4,2)
\psframe[linestyle=solid,linewidth=.5mm](4,0)(5,1)
\psframe[linestyle=solid,linewidth=.5mm](5,-1)(6,0)
\rput(.5,1.5){5}
\rput(2.5,2.5){5}
\rput(3.5,1.5){5}
\rput(4.5,0.5){5}
\rput(5.5,-.5){5}
\psframe[fillstyle=solid,fillcolor=yellow](2,3)(3,4)
\psframe[linestyle=solid,linewidth=.5mm](2,3)(3,4)
\psframe[linestyle=solid,linewidth=.5mm](3,2)(4,3)
\psframe[linestyle=solid,linewidth=.5mm](4,1)(5,2)
\psframe[linestyle=solid,linewidth=.5mm](5,0)(6,1)
\psframe[linestyle=solid,linewidth=.5mm](6,-1)(7,0)
\rput(2.5,3.5){6}
\rput(3.5,2.5){6}
\rput(4.5,1.5){6}
\rput(5.5,0.5){6}
\rput(6.5,-.5){6}
\psframe[fillstyle=solid,fillcolor=yellow](1,3)(2,4)
\psframe[fillstyle=solid,fillcolor=yellow](2,4)(3,5)
\psframe[linestyle=solid,linewidth=.5mm](1,3)(2,4)
\psframe[linestyle=solid,linewidth=.5mm](2,4)(3,5)
\psframe[linestyle=solid,linewidth=.5mm](3,3)(4,4)
\psframe[linestyle=solid,linewidth=.5mm](4,2)(5,3)
\psframe[linestyle=solid,linewidth=.5mm](5,1)(6,2)
\psframe[linestyle=solid,linewidth=.5mm](6,0)(7,1)
\psframe[linestyle=solid,linewidth=.5mm](7,-1)(8,0)
\rput(1.5,3.5){7}
\rput(2.5,4.5){7}
\rput(3.5,3.5){7}
\rput(4.5,2.5){7}
\rput(5.5,1.5){7}
\rput(6.5,0.5){7}
\rput(7.5,-.5){7}
\end{pspicture}
\ec
\vspace{1ex}
\caption{\label{blockProcess}Oriented Digital Boiling Process.
The number in a box is the time this box was added and if the box
is colored, then the box was added stochastically according to Rule 3.}
\end{figure}

This process can be readily
visualized by imagining the growth proceeding by the addition
of unit squares starting with the initial square centered at
$(1/2,-1/2)$.  We denote this initial time by placing a $0$ in
this box.  At time $t=1$ a box is added to the right (centered at
$(3/2,-1/2)$) and with probability $p$ a box is added to the top
of the initial box (centered at $(1/2,1/2)$).  We place $1$'s in
the boxes added at time $t=1$.  The boxes that are added
stochastically (Rule (3)) are shaded.  An example of this process
run for seven time steps is shown in Fig.~\ref{blockProcess}.

\subsubsection{Path Description}

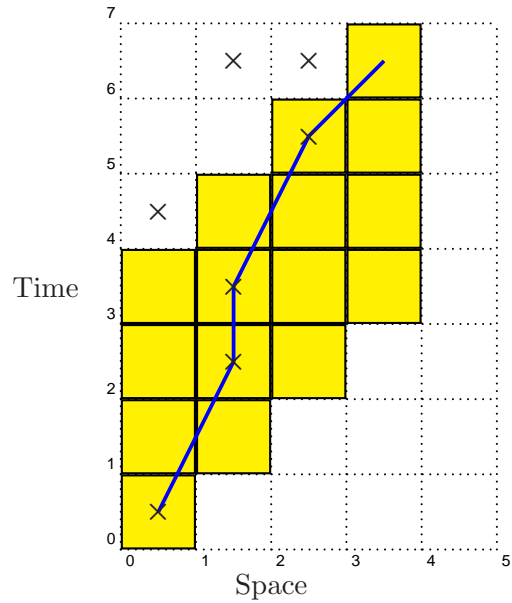
\begin{figure}
\vskip-1cm
\bc
\begin{pspicture}(0,0)(5,7)\showgrid
\psframe[fillstyle=solid,fillcolor=yellow](3,6)(4,7)
\psframe[fillstyle=solid,fillcolor=yellow](3,5)(4,6)
\psframe[fillstyle=solid,fillcolor=yellow](3,4)(4,5)
\psframe[fillstyle=solid,fillcolor=yellow](3,3)(4,4)
\psframe[fillstyle=solid,fillcolor=yellow](2,5)(3,6)
\psframe[fillstyle=solid,fillcolor=yellow](2,4)(3,5)
\psframe[fillstyle=solid,fillcolor=yellow](2,3)(3,4)
\psframe[fillstyle=solid,fillcolor=yellow](2,2)(3,3)
\psframe[fillstyle=solid,fillcolor=yellow](1,4)(2,5)
\psframe[fillstyle=solid,fillcolor=yellow](1,3)(2,4)
\psframe[fillstyle=solid,fillcolor=yellow](1,2)(2,3)
\psframe[fillstyle=solid,fillcolor=yellow](1,1)(2,2)
\psframe[fillstyle=solid,fillcolor=yellow](0,3)(1,4)
\psframe[fillstyle=solid,fillcolor=yellow](0,2)(1,3)
\psframe[fillstyle=solid,fillcolor=yellow](0,1)(1,2)
\psframe[fillstyle=solid,fillcolor=yellow](0,0)(1,1)
\psline[linewidth=0.5mm,
linecolor=blue](.5,.5)(1.5,2.5)(1.5,3.5)(2.5,5.5)(3.5,6.5)
\rput(0.5,.5){\boldmath{$\times$}}
\rput(0.5,4.5){\boldmath{$\times$}}
\rput(1.5,2.5){\boldmath{$\times$}}
\rput(1.5,3.5){\boldmath{$\times$}}
\rput(1.5,6.5){\boldmath{$\times$}}
\rput(2.5,5.5){\boldmath{$\times$}}
\rput(2.5,6.5){\boldmath{$\times$}}
\rput(2.0,-0.5){Space}
\rput(-1.0,3.5){Time}
\end{pspicture}
\ec
\vspace{2ex}
\caption{\label{spaceTime} The backwards lightcone of the point
$(x,t)=(3,7)$
for the process shown in Fig.~\ref{blockProcess}. The $\times$'s
denote the marked points and
polyogonal line gives a longest
increasing path. The length of this path is equal to the
number of $\times$'s in the path.  This length equals $h_t(x)$.}
\end{figure}

As has been observed many times before
(see, e.g., \cite{aldous1, cggk, gravner2, griffeath, praehofer1,
sepp}),
a most productive way to
think about height processes is to introduce the (discrete) backwards
lightcone of a point $(x,t)$.  Precisely, if
$\mathcal{S}=\textbf{Z}_+\times\textbf{Z}_+$ denotes space-time,
then
\[\mathcal{L}_B(x,t) =\{(x',t')\in \mathcal{S}:
0\le x' \le x, x' \le t' < x'+t-x \}.
\]
For those space-time points in $\mathcal{L}_B$
 at which a box was added stochastically
(according to Rule (3)), we place a $\times$.
We call such space-time points \textit{marked}.
 We define the length of a sequence
$\pi=\{(x_1,t_1),\ldots,(x_k,t_k)\}$ of
(distinct) space-time points in $\mathcal{L}_B$  to be $k$.
Such a sequence
  $\pi$ is \textit{increasing}
if $0\le x_i-x_{i-1}\le t_i-t_{i-1}-1$ for $i=2,\ldots,k$.
Let $L(x,t)$ equal the length of the longest increasing
sequence of marked space-time points in $\mathcal{L}_B(x,t)$.
(If $x\le t$ and
$\mathcal{L}_B(x,t)$ contains no increasing
path, then $L(x,t):=0$.)
For the example in
 Fig.~\ref{blockProcess}, the discrete backwards
 lightcone $\mathcal{L}_B(3,7)$ and an
 increasing path are shown in Fig.~\ref{spaceTime}.  One observes
 that $h_7(3)=L(3,7)=4$.  Indeed, this
  is a general fact.  However,
  before proceeding with its proof,
  it is useful to change slightly the point of view of the
 process defined by Rules (1)--(3). Let $\Pi=\Pi(p)$ be a random
 subset of $\mathcal{S}$ to which every point of $\mathcal{S}$ belongs
 with probability $p$.  We \textit{mark} the points of $\mathcal{S}$
 that belong to $\Pi$.  (Accordingly, we
  call these points \textit{marked}.)
  $L(x,t)$ remains the same; namely, the length of the longest
  increasing sequence of marked space-time points in
  $\mathcal{L}_B(x,t)$.  With regard to the process,
  we may intuitively think
  that all the ``coins'' used in Rule (3) are thrown in advance---of
course,
  many of these are ignored as $x$ at time $t$ 
   may become occupied deterministically by Rule (2).  Precisely,
  Rule (3) is replaced with

  \begin{enumerate}
  \item[($\tilde{3}$)] Otherwise, 
  \[ h_{t+1}(x)=\left\{ \begin{array}{rl}
                h_t(x)+1, & \textrm{if} \> (x,t)\in\Pi, \\
                h_t(x),   & \textrm{if} \> (x,t)\not\in\Pi.
       \end{array} \right. \]
  \end{enumerate}
  We are now ready to prove the last passage property

 \textbf{Proposition}~\cite{gravner2}. $h_t(x)=L(x,t)$.

 \textbf{Proof}.
 We first show that our process is \textit{attractive}\footnote{For
 examples of the kind of exotic shapes that can occur from
 cellular automaton rules without this monotonicity property,
 see \cite{gravnerGriffeath}.}
 in the following sense: Let $\Pi$ and $\Pi^\prime$ be two sets of
 marked points such that $\Pi\subset\Pi^\prime$.  Let $h_t$ evolve using
 $\Pi$ and $h_t^\prime$ using $\Pi^\prime$, then $h_t\le h_t^\prime$ for
 all $t$. For if this were \textit{not} true, then, for some
 $t$, $h_s\le h_s^\prime$, $s\le t$, and $h_{t+1}(x)>h_{t+1}^\prime(x)$
 for some $x$. This, of course, implies that $h_t(x)=h_t^\prime(x)$.
  But then $h_{t+1}(x)=h_t(x)+1$ either because of
 Rule (2); in which case, $h_t^\prime(x-1)\ge
h_t(x-1)>h_t(x)=h_t^\prime(x)$,
 so (by Rule (2)) $h_{t+1}^\prime(x)=h_t^\prime(x)+1$; or,
 because $(x,t)\in\Pi\subset\Pi^\prime$, so again
 $h^\prime_{t+1}(x)=h_t^\prime(x)+1$.  This is a contradiction. Thus
 we've established the attractiveness of our process.

 The property of attractiveness immediately implies $h_t(x)\ge L(x,t)$,
 since any increasing path of length $k$ will, without
 the addition of other
 marked points, cause $h_t(x)\ge k$.

 We now show that $h_t(x)\le L(x,t)$.  We will show, by induction
 on $k$ and $t$, that $h_t(x)=k$ implies there exists an
increasing sequence of marked points of length $k$ in $\mathcal{L}_B(x,t)$.
 This is obviously true for either $t=0$ or $k=0$.  (Note that
 $h_t(x)\ge 0$ means that $x\le t$.)  Now assume the claim has been
 demonstrated for all $k'<k$ and $t'<t$.  We can clearly assume that
 $h_{t-1}(x)=k-1$, or else we can use the induction hypothesis right
 away.  Therefore, we have two possibilities.

 \textit{Case 1.} $h_t(x)=h_{t-1}(x)+1$ by application of Rule (2).
 This means (by Rule (2)) that $h_{t-1}(x-1)=k$.  Thus by the
 induction hypothesis, there is an increasing sequence of
 length $k$
 in $\mathcal{L}_B(x-1,t-1)\subset \mathcal{L}_B(x,t)$.

 \textit{Case 2.}  $h_t(x)=h_{t-1}(x)+1$ by application of Rule
($3^\prime$).
 This means that $h_{t-1}(x)=k-1$ and $(x,t-1)\in\Pi$.
   By the induction
 hypothesis, $\mathcal{L}_B(x,t-1)$ contains an increasing
 path of length $k-1$.  Adjoin the marked point $(x,t-1)$ to
 the sequence.  Observe that the increasing property is
 preserved.  This completes the proof of the proposition.

We summarize this section by noting that $h_t$ satisfies
for all $t\ge 1$, $x\ge 0$, 
\[ h_t(x) =\max\left\{h_{t-1}(x-1),h_{t-1}(x)+\epsilon_{x,t}\right\} \]
where $\epsilon_{x,t}=1$ if $(x,t)\in\Pi$ and $0$ otherwise.
The initial conditions are (\ref{initial}).  (We take $h_t(-1)=-\iy$.)
Formulated this way ODB is a ``stochastic dynamic programming'' problem.
\subsubsection{The $(0,1)$-Matrix Description of ODB}

Without changing the increasing path property, the
backwards lightcone $\mathcal{L}_B$ of any space-time
point $(x,t)$ can be deformed into a rectangle of
size $(t-x)\times (x+1)$. Thus the equivalent
problem is to fix $x$ and $t$ and to set
$m=t-x$,  $n=x+1$,  and to
consider a $(0,1)$-matrix $A$ of size $m\times n$. We number
the rows of $A$ starting at the bottom
of $A$ and the columns of $A$ starting
at the left of $A$. A increasing path in $\mathcal{L}_B$
becomes
a sequence of $1$'s in $A$ at, say,  positions
$\{(i_1,j_1),\ldots,(i_k,j_k)\}$
such that the $i_\ell$ ($\ell=1,\ldots,k$) are
increasing and the $j_\ell$ ($\ell=1,\ldots,k$) are weakly
increasing.
Any such $(0,1)$-matrix $A$ of size $m\times n$ corresponds
(bijectively) to a two-line array\footnote{We have chosen both a
nonstandard labeling of $A$ and a
nonstandard bijection $A\leftrightarrow w_A$
so that our increasing path property remains (essentially) the same
under the bijections.}
\be w_A=\left(
\begin{array}{cccc}
j_1&j_2&\cdots&j_k\\
i_1&i_2&\cdots&i_k
\end{array}
\right)\label{twoLineArray}\ee
where $j_1\le j_2 \le \cdots\le j_k$ and if $j_\ell=j_{\ell+1}$, then
$i_{\ell} < i_{\ell+1}$ and the pair ${j\choose i}$ appears in $w_A$
if and only if the $(i,j)$ entry of $A$ is 1.
Note that the upper numbers belong to
$\{1,2,\ldots,n\}$ and the lower numbers to
$\{1,2,\ldots,m\}$.
For example, the matrix
\[
A=\left(\begin{array}{ccccccc}
        0&0&0&1&0&0&1\\
        1&1&1&1&0&1&1\\
        1&1&0&0&0&1&0\\
        1&0&1&1&0&1&1\\
        0&1&1&0&0&0&0\\
        0&0&0&0&1&0&1
        \end{array}\right)
\]
maps to the two-line array
\begin{eqnarray*} w_A&=&\left(
\begin{array}{cccccccccccccccccccc}
1&1&1&\textbf{2}&2&2&3&
\textbf{3}&3&4&4&4&5&6&\textbf{6}&\textbf{6}&7&7&7&\textbf{7}
\\
3&4&5&\textbf{2}&4&5&2&\textbf{3}&5&3&5&6&1&3&
\textbf{4}&\textbf{5}&1&3&5&\textbf{6}
\end{array}\right).
\end{eqnarray*}
(Recall the convention for row labels.)
As an example, a
 longest increasing path (of length 5) is indicated
in bold typeface.  We remark that one can compute the
length of an increasing path by \textit{patience
sorting}~\cite{aldous2}  on
the bottom row  of $w_A$ (from left to right) with the rule that a
number is placed
on the left most pile such that it is less than or equal
to the number showing in the pile.  Patience sorting on
the above example results in the five piles
\[\begin{array}{ccccc}
&3&&&\\
1&3&4&&\\
1&3&5&&\\
2&3&5&5&\\
2&4&5&5&\\
3&4&5&6&6.
\end{array}\]

If $N$ denotes the number of $1$'s in a random $m\times n$
$(0,1)$-matrix
$A$, then the above mappings imply that for any nonnegative integer $h$,
\bae
\pr\left(h_t(x)\le h\right)&=&\sum_{k\ge 0} \pr\left(h_t(x)\le h\vert
N=k\right) \pr\left(N=k\right)\nonumber\\
&=&\sum_{k=0}^{mn} {mn\choose k} p^k (1-p)^{mn-k} \,
\pr\left(L_{m,n,k}\le h\right)\label{distr1}
\eae
where $L_{m,n,k}$ is the length of the longest increasing path
in a random (0,1)-matrix $A$
with $k$ 1's (or equivalently, in the associated $w_A$).

\subsubsection{Tableaux Description of ODB}
The dual RSK algorithm~\cite{knuth, stanley}
 is a bijection between $(0,1)$-matrices
$A$ of size $m\times n$ and pairs $(P,Q)$ such that
$P^t$ (the transpose of $P$) and $Q$ are 
semistandard Young tableaux (SSYTs) with
$\textrm{sh}(P)=\textrm{sh}(Q)$ where the elements of $P$ are from
$\{1,2,\ldots,m\}$ and the elements of $Q$ are
from $\{1,2,\ldots,n\}$.  In terms of the associated
$w_A$, (\ref{twoLineArray}),  one forms $P$ by successive row bumping
of the second row of $w_A$ starting with $i_1$ and with the rule
an element $i$ bumps the leftmost element $\ge i$.  Thus each
row of $P$ is strictly increasing.  A fundamental property of
the dual RSK algorithm is that the
length of the longest strictly increasing subsequence of the
second row of $w_A$ equals the number of boxes in the first row of $P$.

If $d_\la(M)$ denotes the number of SSTYs of shape $\la$ with entries
coming from $\{1,2,\ldots,M\}$, then the number of pairs $(P,Q)$
of fixed shape $\la$ in the above dual RSK algorithm is
\[ d_{\la'}(m) d_\la(n) \]
where $\la'$ is the conjugate partition.  (Conjugate since $P^t$ is
a SSYT.)  Since there are ${mn\choose k}$ $(0,1)$-matrices
with $k$ 1's,
\[
\pr\left(L_{m,n,k}\le h\right) ={1\ov {mn\choose k}}\sum_{{\la\vdash k
\atop \la_1\le h}}
d_{\la'}(m) d_\la(n)= {1\ov {mn\choose k}}\sum_{{\la\vdash k \atop
\ell(\la)\le h}} d_\la(m) d_{\la'}(n). \]
And hence from (\ref{distr1})
\[
\pr\left(h_t(x)\le h\right)=(1-p)^{mn}\,\sum_{k=0}^{mn} r^k
\sum_{{\la\vdash k \atop \ell(\la)\le h}} d_\la(m) d_{\la'}(n)\]
where $r=p/(1-p)$.  Observe that for $\vert\la\vert>mn$,
$d_\la(m)d_{\la'}(n)=0$.
(A SSYT with entries from $\{1,2,\ldots,M\}$ can have at most $M$ rows.)
If $\mathcal{P}$ denotes the set of all partitions (including
the empty partition), then the
 above sum can be summed over all partitions without changing its value,
 \be
 \pr\left(h_t(x)\le h\right) = (1-p)^{mn}
 \sum_{{\la\in\mathcal{P}\atop \ell(\la)\le h}} r^{|\la|}d_\la(m)
d_{\la'}(n).
  \label{distr2}\ee

Comparing (\ref{distr2}) with Johansson's Krawtchouck ensemble results
establishes the equivalence of ODB with the 
Sepp{\"a}l{\"a}inen-Johansson model.

\subsubsection[Gessel's Theorem and the Borodin-Okounkov
Identity]{Application of Gessel's Theorem
and the Borodin-Okounkov Identity}

 Gessel's theorem \cite{gessel,tw5} is
 \[ \sum_{{\lambda\in\mathcal{P}\atop \ell(\lambda)\le h}}
  r^{\vert\lambda\vert}
 s_{\lambda}(x) s_{\lambda}(y) = D_h(\vp) \]
 where $s_\la$ are the Schur functions (see, e.g.~\cite{stanley}) and
 $D_h(\vp)$ is the $h\times h$ Toeplitz
 determinant\footnote{If $\phi$ is a function on the unit circle
 with Fourier coefficients $\phi_k$ then $T_n(\phi)$ denotes
 the Toeplitz matrix $(\phi_{i-j})_{i,j=0,\ldots,n-1}$ and
 $D_n(\phi)$ its determinant.}
 with symbol
 \[ \vp(z) = \prod_{j=1}^\infty (1-x_j  z)^{-1} \prod_{j=1}^\infty
 (1-y_j r z^{-1})^{-1} .\]
 If we apply to
 both sides of this identity the automorphism $\omega$ (see
  Stanley~\cite{stanley}, pg.~332),
 $\omega(s_\lambda)=s_{\lambda^\prime}$, to the symmetric functions in
 the $x$-variables we obtain
 \be \sum_{{\lambda\in\mathcal{P}\atop \ell(\lambda)\le h}}
  r^{\vert\lambda\vert}
 s_{\la^\prime}(x) s_{\la}(y) = D_h(\vp)
 \label{gessel}\ee
 where \textit{now} the symbol is
 \be \varphi(z) = \prod_{j=1}^\infty(1+x_j  z) \prod_{j=1}^\infty
 (1-y_j r z^{-1})^{-1}.\label{symbol1}\ee

 Recalling the specialization
 $\textrm{ps}_n^1$ (see Stanley~\cite{stanley}, pg.~303), we apply
 $\textrm{ps}_n^1$ to the $x$-variables and $\textrm{ps}_m^1$ to
 the $y$-variables in
 Gessel's identity (\ref{gessel}) and observe\footnote{Note that
 $\textrm{ps}_n^1 \,s_\la=d_\la(n)$ which follows
 from the combinatorial definition of the Schur function.}
 that the resulting LHS is precisely the RHS
 of (\ref{distr2}).  Since the specialization $\textrm{ps}_n^1$
 is a ring homomorphism, we
 may apply it directly to the
 symbol (\ref{symbol1}).  Doing so we  obtain

\be \pr\left(h_t(x)\le h\right)=(1-p)^{mn}
D_h(\vp) \label{distr3}\ee
where
\be \vp(z)=(1+z)^n (1-r/z)^{-m}. \label{symbol}\ee
This derivation required $r<1$.  However, by (\ref{distr1}) the left
side is a rational function of $r$, and analytic continuation
shows that (\ref{distr3}) holds for all $r\ge0$ if in the integral
representing the Fourier coefficients of $\vp$ the contour has
$r$ on the inside.

The Borodin-Okounkov~\cite{borodin2}
 identity expresses a Toeplitz determinant
in terms of a Fredholm determinant of an infinite matrix which
in turn is a product of two Hankel matrices.  Subsequent
simplifications of the proof by Basor and Widom \cite{basorWidom}
extended the identity to block Toeplitz determinants.  We now
apply this identity to the Toeplitz determinant (\ref{distr3}).
First we find the Wiener-Hopf factorization of $\varphi(z)$:
\[ \vp(z)=\vp_+(z)\,\vp_-(z) \]
where
\[ \vp_+(z)=(1+z)^n\, , \>\>\> \vp_-(z)=(1-r/z)^{-m}.\]
Define $K_h$ acting on $\ell^2(\{0,\,1,\cd\})$ by
\be
K_h(j,k)=\sum_{\l=0}^{\iy}(\vp_-/\vp_+)_{h+j+\l+1}\;
(\vp_+/\vp_-)_{-h-k-\l-1}.
\label{Kdefn}\ee
The Borodin-Okounkov identity is then
\[ D_h(\vp)= Z\,\det\left(I-K_h\right). \]
Since the determinant on the right tends to 1 as $h\ra\iy$ as
does $\pr\left(h_t(x)\le h\right)$, we have $Z=(1-p)^{-mn}$.
Thus we have derived a representation of the distribution function
of the random variable $h_t(x)$
in terms of a Fredholm determinant,
\be
\pr\left( h_t(x)\le h \right)=\det\left( I-K_h
\right).\label{distrFred}\ee
This derivation also required $r<1$.  As above, analytic continuation
shows that (\ref{distrFred}) holds for all $r\ge0$ if in the integral
representing the Fourier coefficients of $\vp_{-}/\vp_{+}$ the contour
has $r$ on the inside and $-1$ on the outside.  (In fact the contour
must have $-1$ on the outside no matter what $r$ is.)

A somewhat different direction (and one we do not follow here) is
to apply  isomonodromy and Riemann-Hilbert methods~\cite{bdj, deift,
itw2}
 directly to the
Toeplitz determinant $D_h(\vp)$.  This would result in the
identification of
$D_h(\vp)$ as a $\tau$-function of an integrable ODE.

\subsection{Inhomogeneous ODB}
In ODB the probability $p$ appearing in Rule (3) is independent of
the site $x$.  \textit{Inhomogeneous ODB} replaces
Rule (3), for each site
$x\in\textbf{Z}_+$,  with
\par
\par
$(3_x)$ Otherwise,  then independently
of the other sites and other times,
$h_{t+1}(x)=h_t(x)+1$ with probability $0<p_x<1$ and $h_{t+1}(x)=h_t(x)$
with probability $q_x:=1-p_x$.
\par
\par
Since the dual RSK algorithm
is a bijection between $(0,1)$-matrices $A$ and pairs $(P,Q)$ such
that $P^t$ and $Q$ are SSYTs \textit{with} $\textrm{col}(A)=\textrm{type}(P)$
and $\textrm{row}(A)=\textrm{type}(Q)$ \cite{stanley}, we have
\be
\pr\left(h_t(x)\le h\right)=q_0^m \cdots q_x^m\,
\sum_{{\la\in\mathcal{P}\atop \ell(\la)\le h}}d_{\la}(m)\, s_{\la^\prime}(r)
\label{inhomoMeasure}\ee
where, as before, $m=t-x$, but now
$ r=\left(r_0,\ldots,r_x,0,\ldots\right)$
with $r_j:=p_j/q_j$.
The proof of (\ref{inhomoMeasure}) is
straightforward and
similar to the proof of the analogous
result in \cite{itw1}; therefore, we omit it. The right hand
side of (\ref{inhomoMeasure}) clearly  reduces
to (\ref{distr2}) in the homogeneous case.

We  again apply Gessel's theorem to obtain the Toeplitz
determinant representation
\[ \pr\left(h_t(x)\le h\right)=q_1^m \cdots q_n^m\, D_h(\vp) \]
where\footnote{The homogeneous case of (\ref{inhomoSymbol}) does
not directly reduce to (\ref{symbol}).  It does after $z\ra z/r$ which
corresponds to a similarity transformation of the Toeplitz matrix.}
\be \vp(z) = (1-1/z)^{-m}\, \prod_{j=0}^x (1+r_j z)\, .
 \label{inhomoSymbol}\ee
Application of the Borodin-Okounkov identity  results in
 a Fredholm determinant representation for this distribution function.
Observe that from either (\ref{inhomoMeasure})  or (\ref{inhomoSymbol})
it follows that \textit{$\pr(h_t(x)\le h)$ is a symmetric function
of $(p_0,p_1,\ldots,p_x)$}.  This property opens the possibility
for an analysis of the spin glass version of ODB which we plan to
address in future work.

\subsection{Weak ODB and Strict ODB}

Here are two natural variants of the ODB.  We let
the ``spontaneous increase'' in Rule (3) apply after Rule (2) has
already
taken effect to get \textit{weak ODB}:
\begin{enumerate}
\item[($1'$)] $h_t(x)\le h_{t+1}(x)$ for all space-time points $(x,t)$.
\item[($2'$)] If $h_t(x-1)>h_t(x)$, then $\tilde{h}_t(x)=h_t(x-1)$ else
$\tilde{h}_t(x)=h_{t-1}(x)$.  (Here $\tilde{h}_t$ is an intermediate
height function.)
\item[($3'$)] Independently of the other sites and other times,
$h_{t+1}(x)=\tilde{h}_t(x)+1$ with probability $p$.  (With probability
$1-p$, $h_{t+1}(x)=\tilde{h}_t(x)$.)
\end{enumerate}

In \textit{strict ODB} we require that the left neighbor is 
\textit{rested}\footnote{The height at a site cannot increase
at two consecutive times, i.e.\ it must rest for one time unit
before it is allowed to increase.} for
the spontaneous increase.  
(We take $h_t(x)=-\iy$ for $x<0$ which in this model implies
$h_t(0)=0$ for every $t$.)
\begin{enumerate}
\item[($1''$)] $h_t(x)\le h_{t+1}(x)$ for all space-time points $(x,t)$.
\item[($2''$)] If $h_t(x-1)>h_t(x)$, then $h_{t+1}(x)=h_t(x-1)$.
\item[($3''$)] Otherwise, if $x-1$ is rested at time $t$,
$h_t(x-1)=h_t(x)$ then
independently of other
sites and times, $h_{t+1}(x)=h_t(x)+1$ with probability $p$
($h_{t+1}(x)=
h_t(x)$ with probability $1-p$.)
\end{enumerate}

In a similar way one shows
\begin{itemize}
\item In weak ODB, $h_t(x)$ equals, in distribution, the
longest sequence $(i_\ell,j_\ell)$ of
positions in a random $(0,1)$-matrix of size $m\times n$
($m=t-x+1$, $n=x+1$) which have
entry $1$ such that  $i_\ell$ are $j_\ell$ are both weakly increasing.
(The lower left corner of the matrix is fixed to be a 0.)
\item In strict ODB, $h_t(x)$ equals, in distribution, the
 longest sequence $(i_\ell,j_\ell)$ of
positions in a random $(0,1)$-matrix of size $m\times n$ ($m=t-x$,
$n=x$) which have
entry $1$ such that  $i_\ell$ are $j_\ell$ are both strictly
 increasing.
\end{itemize}

\section{Limit Theorems}
\setcounter{equation}{0}
\begin{figure}
\bc
\hspace{1cm}\resizebox{12cm}{5cm}{\includegraphics{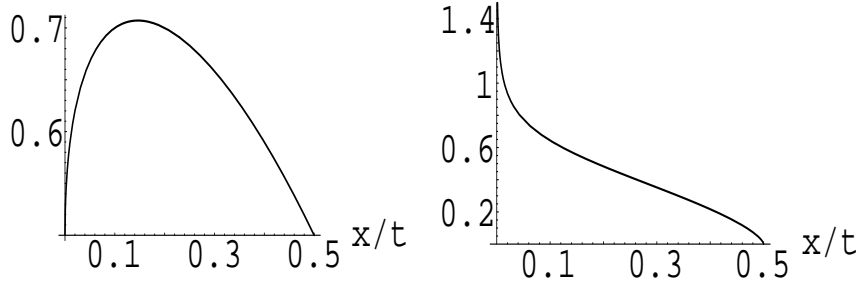}}
\caption{\label{c1c2Fig} In the GUE Universal Regime, the left figure
displays the limiting shape, $c_1$,
as a function of
$x/t$ and the right figure displays the normalization constant,
$c_2$, as a function of $x/t$. In both cases, $p=1/2$.}
\ec
\end{figure}

\begin{figure}
\bc
\hspace{2cm}\resizebox{9cm}{6cm}{\includegraphics{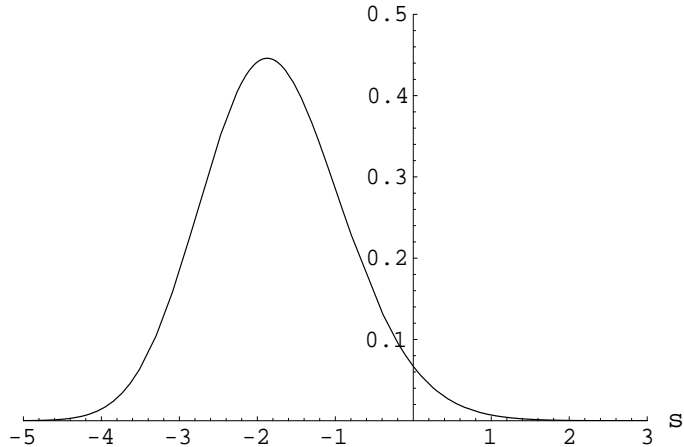}}
\caption{\label{TWdistr}The density $f_2(s)=dF_2/ds$ where
$F_2$ is defined by (\ref{F2}). The distribution function $F_2$
has mean $\mu=-1.77109$, standard deviation $\s=0.9018$, skewness
$S=0.2241$
and excess kurtosis $K=0.0935$.}
\ec
\end{figure}

In this section we derive limit theorems
for the distribution function
 $\pr\left(h_t(x)\le h\right)$ for ODB.  Our
 starting point will be the Fredholm determinant representation
 (\ref{distrFred}).   This distribution function is a function
 of four variables, $x$, $t$, $h$
 and $p$; and accordingly, there
 are several asymptotic regimes:

 \begin{enumerate}
\item[(1)]  \textit{GUE Universal Regime}:
Let $x\ra\iy$, $t\ra\iy$ such that $p_c:=1-x/t<1$ is fixed.  For
fixed $p<p_c$ define
\be c_1:=2p_c p - p +2\sqrt{p p_c (1-p)(1-p_c)}, \label{c1}\ee
\be c_2:=(p_c(1-p_c))^{1/6} (p(1-p))^{1/2}\left[
\left(1+\sqrt{{(1-p)(1-p_c)\ov p p_c}}\right)\left(
\sqrt{{p_c\ov 1-p_c}}
-\sqrt{{p\ov 1-p}}\right)\right]^{2/3}\, .\label{c2}\ee
 We will show that
\[ \pr\left({h_t(x)-c_1 t\ov c_2\, t^{1/3}}< s\right)
\ra F_2(s) \]
where~\cite{tw1}
\be F_2(s)=\det\left(I-K_{\textrm{Airy}}\right)=\exp\left(-\int_s^\iy
(x-s) q(x)^2\,dx\right).
\label{F2}\ee
Here $K_{\textrm{Airy}}$ is the operator with Airy kernel
acting on $L^2((s,\iy))$ (see (\ref{airyKernel}) below)
and $q$ is the (unique) solution of the Painlev\'e II equation
\[ q''=sq+2q^3 \]
with boundary condition $q(s)\sim\textrm{Ai}(s)$ as $s\ra\iy$.  The
limiting shape, $c_1$, and
the normalization constant, $c_2$,  as  functions of $x/t$ are shown in
Fig.~\ref{c1c2Fig} for $p=1/2$.
The probability
density, $f_2=dF_2/ds$,
  is shown in Fig.~\ref{TWdistr}.
\item[(2)] \textit{Critical Regime}: Let $x\ra\iy$, $t\ra\iy$ such
that
\[ x=(1-p)t+o(\sqrt{t}). \]
For fixed $\Delta\in\textbf{Z}_+$ we will show    that
\[ \pr\left(h_t(x)-(t-x)\le -\Delta \right) \]
converges to a $\Delta  \times \Delta $ determinant.
One can think of this as
\[ p=p_c+o({1\ov\sqrt{t}}).\]
\item[(3)] \textit{Deterministic Regime}: For $x\ra\iy$,
$t\ra\iy$ and fixed $p>p_c$,   we will
show that
\[ \pr\left(h_t(x)=p_c t\right)\ra 1.\]
\item[(4)] \textit{Finite $x$ GUE Regime}: Fix $x$ and let $t\ra\iy$,
then we will show that
\[ \pr\left({h_t(x)-p\, t\ov (p(1-p)\, t)^{1/2}}< s\right) \]
converges to the distribution of the largest eigenvalue in the
GUE of $(x+1)\times(x+1)$ hermitian matrices, denoted below by
$F_{x+1}^{GUE}$.
\end{enumerate}

\subsection{GUE Universal Regime}
It is convenient to use the variables $m=t-x$ and $n=x+1$ rather
than $x$ and $t$ and to translate back to the space-time variables
at the end.  We assume $p<p_c:=m/(n+m)$.  (This is asymptotically
$1-x/t$ as defined above.)
Further, when there is no chance of confusion, we denote the
random variable $h_t(x)$ by $H$.  (We reserve lower case
 $h$  to denote the values of $H$.)
Set $h=c\,m+s\,m^{1/3}$, where $c$ will be determined shortly, and
$\al=n/m$.  (In this notation the condition $p<p_c$ is $\al\, r<1$.)
For any $v$ the matrix $((-v)^{k-j}\,K_h(j,k))$ has the same Fredholm
determinant
(the determinant of $I$ minus the matrix) as $(K_h(j,k))$. We shall show
that for
a particular $v$ and a certain constant $b>0$ this matrix scales to a
kernel with
the same Fredholm determinant as
\be K_{{\rm Airy}}(s/v(3b)^{1/3}+x,\,s/v(3b)^{1/3}+y),\label{KAiry}\ee
on $(0,\,\iy)$, where
\be K_{{\rm Airy}}(s+x,\,s+y)=\int_0^{\iy}{\rm Ai}(t+s+x)\,{\rm
Ai}(t+s+y)\,dt.\label{airyKernel}\ee
This gives
\be\lim_{m\ra\iy}\,{\rm Prob}\left({H-c\, m\ov m^{1/3}}\le
s\right)=F_2(s/v(3b)^{1/3}).\label{problim}\ee

Here is what we mean by scaling. Any matrix $(M(j,k))$ acting on
$\l^2(\textbf{Z}_+)$
has the same Fredholm determinant as the kernel $M([x],\,[y])$ on
$L^2(0,\,\iy)$ and
this in turn has the same Fredholm determinant as
$M_m(x,\,y)=m^{1/3}\,M([m^{1/3}\,x],\,[m^{1/3}\,y])$.
If this kernel has the limit $k(x,\,y)$ we say that the matrix
$(M(j,k))$ has, after
the scaling $j\ra m^{1/3}\,x,\ \ k\ra m^{1/3}\,y$, the limit $k(x,\,y)$.
If $M_m(x,\,y)$
converges to $k(x,\,y)$ in trace norm then the Fredholm determinant of
$(M(j,k))$
converges to that of $k(x,\,y)$. And if $(M(j,k))$ were the product of
two matrices
each having scaling limits in Hilbert-Schmidt norm (under the same
scaling, of
course), then
the Fredholm determinant of the product converges to the Fredholm
determinant of the
product of the limits. This is what we shall show in our case.

There is a slightly awkward notational problem. Since $h$ is always an
integer
and $h=c\,m+s\,m^{1/3}$, the quantity $s$ as it appears here and the
analysis which
follows is not completely arbitrary. What we actually show is
that if $h$ and $m$ tend to infinity, and $s$ is defined in terms of
them
by the formula $h=c\,m+s\,m^{1/3}$, then
\be{\rm Prob}\left(H\le
h\right)-F_2\left(s/v(3b)^{1/3}\right)\ra 0\label{probdif}\ee
uniformly for $s$ lying in a bounded set. From this we easily deduce
(\ref{problim})
for fixed $s$,
which now has a different meaning. These observations are important when
one tries to
estimate errors. It can be shown that the difference in (\ref{probdif})
is $O(m^{-2/3})$.
But the
difference between the right side of (\ref{problim}) and the probability
on the
left can only be expected to be $O(m^{-1/3})$. The reason is that
if the quantity $s'$ is defined by $cm+s'm^{1/3}=[cm+sm^{1/3}]$ then the
probability
is within $O(m^{-2/3})$ of $F_2(s'/v(3b)^{1/3})$, but $s-s'$ is very
likely of the order $m^{-1/3}$.

\subsubsection{The Saddle Point Method}
The matrix $(K_h(j,k))$ is the
product of two matrices, the matrix on the right having $j,k$ entry
$(\ph_+/\ph_-)_{-h-j-k-1}$
and the one on the left having $j,k$ entry
$(\ph_-/\ph_+)_{h+j+k+1}$.
Notice that the first vanishes if $h+j+k+1>m$ so we may assume that all
our indices
$j$ and $k$ satisfy $h+j+k<m$.
We have
\bae (\ph_+/\ph_-)_{-h-j-k-1}&=&{1\ov 2\pi
i}\int(1+z)^n\,(z-r)^m\,z^{-m+h+j+k}\,dz \label{phi+-}\\
&=&(-1)^{h+j+k}{1\ov 2\pi i}\int(1+z)^n\,(r-z)^m\,(-z)^{-m+h+j+k}\,dz,
\nonumber\eae
and a similar formula holds for $(\ph_-/\ph_+)_{h+j+k+1}$. If we set
\[\ps(z)=(1+z)^n\,(r-z)^m\,(-z)^{-(1-c)\,m}\]
then
\[(-1)^{h+j+k}(\ph_+/\ph_-)_{-h-j-k-1}={1\ov 2\pi i}\int
\ps(z)\,(-z)^{s\,m^{1/3}+j+k}\,dz\]
and
\[(-1)^{h+j+k}(\ph_-/\ph_+)_{h+j+k+1}={1\ov 2\pi i}\int
\ps(z)\inv\,(-z)^{-s\,m^{1/3}-j-k-2}\,dz.\]

The contours for the first integral surrounds 0 while the contour for
the
second integral has $r$ on the inside and $-1$ on the outside. The
restriction
$h+j+k<m$ is the same as $s\,m^{1/3}+j+k<(1-c)\,m$.
If we make the replacements $j\ra m^{1/3}\,x,\ \ k\ra m^{1/3}\,y$ these
become
\[{1\ov 2\pi i}\int \ps(z)\,(-z)^{m^{1/3}\,(s+x+y)}\,dz,\ \ \
{1\ov 2\pi i}\int
\ps(z)\inv\,(-z)^{-m^{1/3}\,(s+x+y+2\,m^{-1/3})}\,dz.\]
Our restrictions become $s+x+y<(1-c)\,m^{2/3}$. For convenience we
replace
$s+x+y$ by $x$, and we want to evaluate
\be{1\ov 2\pi i}\int \ps(z)\,(-z)^{m^{1/3}\,x}\,dz,\ \ \
{1\ov 2\pi i}\int \ps(z)\inv\,(-z)^{-m^{1/3}\,x-2}\,dz\label{psints}\ee
asymptotically.  Our restriction is now
$x<(1-c)\,m^{2/3}$.

\begin{figure}
\bc
\resizebox{8cm}{8cm}{\includegraphics{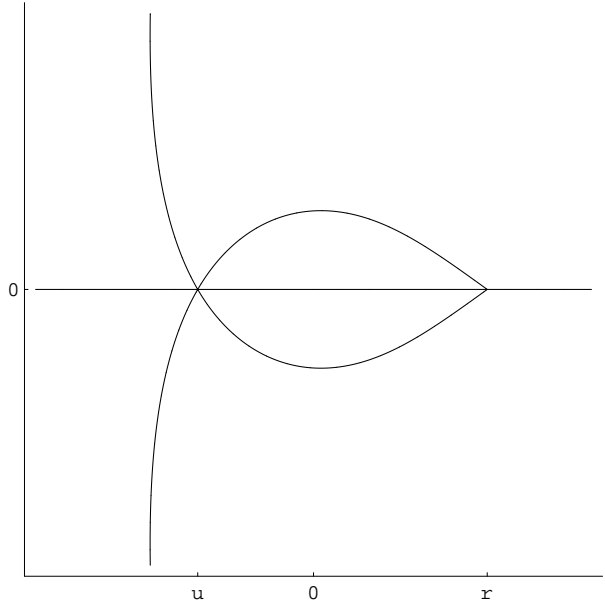}}
\caption{\label{sdcurves}  The steepest descent curves $C^{\pm}$ as described in
the text.}
\ec
\end{figure}

To do a steepest descent we have to find the zeros of
\[
{d\ov dz}\log\,\ps(z)={n\ov 1+z}+{m\ov z-r}-{(1-c)\,m\ov z},\]
or equivalently  the zeros of
\[(c+\al)\,z^2+(c+r-c\,r - \al\,r)\,z+r\,(1-c).\]
(Recall that $\al=n/m$.)
The discriminant of this quadratic
equals zero when
\be c={1\ov 1+r}\,\left(2\sqrt{\al\,r}\,
+(1-\al)\,r\right).\label{velocity}\ee
This is the value of $c$ we take.\footnote{If there were
two critical points, or if we took the negative square root in
(\ref{velocity}), the Fredholm determinant would
tend exponentially to either zero or one. It is only for this value of
$c$ that
we get a nontrivial limit.}
The critical probability is the condition $c=1$, i.e.\ $p_c=m/(m+n)$.
The single zero of the quadratic
is then at $u=-v$ where
\[v={(1-r)\,c + (1-\al)\,r\ov 2\,(c+\al)}={1-\sqrt{\al\, r}
\ov 1+\sqrt{\al/r}}.\]
Note that $0<v<1$ since $0<p<p_c$.
(It is because $u<0$ that we used powers of $-z$
rather than $z$
in the definition of $\ps$.) We write
\ba
6 b:={1\ov m}{d^3\ov
dz^3}\log\,\ps(z)\Big|_{z=u}&=&{2\al\ov(1+u)^3}+{2\ov (u-r)^3}-
{2(1-c)\ov  u^3}\\
&=& {2\left(\sqrt{\al}\,+\,\sqrt{r}\right)^5\ov r\,
\sqrt{\al}\, (1+r)^3 \, (1-\sqrt{\al r}) }.
\ea
The quantity $b$ is positive since $\al\, r<1$.
In the neighborhood
of $z=u$,
\be\ps(z)\sim\ps(u)\,e^{m b \, (z-u)^3}.\label{psasym}\ee
The steepest descent curves will come into $u$ at angles $\pm\pi/3$
and $\pm2\pi/3$. Call the former $C^+$ and the latter $C^-$.
For the integral involving $\ps(z)$
we want $|\ps(z)|$ to have a maximum at that point of the curve and for
the integral involving $\ps(z)\inv$ we want $|\ps(z)|$ to have a minimim
there.
Since $b>0$ the curve for $\ps(z)$ must be $C^+$ and the curve for
$\ps(z)\inv$
must be $C^-$. Both contours will be
described downward near $u$. The curve $C^+$ will loop
around the origin and close at $r$, the upper and lower parts making an
angle there depending on $c$ while $C^-$ will loop
around on both sides and go to infinity with slopes depending
on $c$. (That $C^{\pm}$ have these forms follows from the fact
that the contours cannot cross and, since the only critical
point is at $z=u$, the contours can end only where $\psi$,
respectively $\psi^{-1}$, is zero.)
The steepest descent curves are shown in Fig.~\ref{sdcurves}.

Proceeding formally now, consider the $\ps(z)$ integral and make the
substitution
$z\ra u+z=-v+z$. Then the old $-z$ becomes the new
$v\,(1-z/v)\sim v\,e^{-z/v}$, and recall (\ref{psasym}).
If we make these replacements in the integral we get
\[\ps(u)\,v^{m^{1/3}x}{1\ov 2\pi i}\int_{\iy e^{i\pi/3}}^{\iy
e^{-i\pi/3}}
e^{mbz^3-m^{1/3}xz/v}\,dz.\]
The contour can be deformed to the imaginary axis since we only
pass through regions
where $\Re z^3$ is negative. If we then set $z=-i\z/m^{1/3}$
the
above becomes
\[-\ps(u)\, m^{-1/3}\,v^{m^{1/3}x}{1\ov
2\pi}\int_{-\iy}^{\iy}e^{i\,b\z^3+i\,x\z/v}\,d\z
=-{\ps(u)\, m^{-1/3}\,v^{m^{1/3}x}\ov (3b)^{1/3}}{\rm
Ai}(x/v(3b)^{1/3}).\]
If we recall that $x$ was a replacement for $s+x+y$ we see that the
matrix with
$j,k$ entry
\[(-1)^{h+j+k}\ps(u)\inv v^{-m^{1/3}s-j-k}(\ph_+/\ph_-)_{-h-j-k-1}\]
has the scaling limit
\[-{1\ov (3b)^{1/3}}{\rm Ai}((s+x+y)/v(3b)^{1/3}).\]
Similarly the matrix with $j,k$ entry
\[(-1)^{h+j+k} \ps(u) v^{m^{1/3}\,s+j+k}(\ph_-/\ph_+)_{h+j+k+1}\]
has  $1/v^2$ times exactly same scaling limit. Hence the scaling limit
of
$u^{k-j}\,K_h(j,k))$, which has the same Fredholm determinant as
$(K_h(j,k))$,
is the product of these scaling limits,
\[{1\ov v^2(3b)^{2/3}}\int_0^{\iy}{\rm Ai}((s+x+t)/v(3b)^{1/3})\,
{\rm Ai}((s+t+y)/v(3b)^{1/3})\,dt\]
\[={1\ov v(3b)^{1/3}}\,K_{{\rm
Airy}}((s+x)/v(3b)^{1/3},\,(s+y)/v(3b)^{1/3}).\]
And, as promised,\footnote{The time constant $c_1=p_c \,c$ and
the normalization constant $c_2=p_c^{1/3}\, v\, (3b)^{1/3}$.  A
computation
then gives (\ref{c1}) and (\ref{c2}).}
 this kernel has the same Fredholm determinant as
\[K_{{\rm Airy}}(s/v(3b)^{1/3}+x,\,s/v(3b)^{1/3}+y).\]

\subsubsection{Convergence Proof}
Now for the justification. We have to obtain not only the pointwise
limit, but
uniform estimates to establish convergence of the operators in trace
norm.
We first obtain asymptotics under the assumption that
$x$ lies in a bounded set. (Notice that $x\ge s$ always.)
We begin with
\be{1\ov 2\pi i}\int_{C^+} \ps(z)\,z^{m^{1/3}\,x}\,dz,\label{psint}\ee
and denote by $C^+_\ep$ the portion of $C^+$ which lies within $\ep$ of
the critical point
$\u$.\sp

\noi{\bf Lemma 1}. If in (\ref{psint}) we integrate only over
$C^+_\ep$ the error incurred
is $O(|\ps(u)|\,e^{-\dl\,m})$ for some $\dl>0$. \sp

\noi{\bf Proof}. Define
\[\si(z)={1\ov
m}\log\,\ps(z)=\al\,\log\,(1+z)+\log\,(r-z)+(c-1)\log\,(-z).\]
Its maximum on $C^+$ (it is real-valued there) is $\si(u)$ and it is
stricly
less than this on the complement of $C_\ep$ in $C^+$. Therefore
$\ps(z)/\ps(u)=O(e^{-\dl\,m})$ for some $\dl>0$ on the complement while
$z^{m^{1/3}x}
=e^{O(m^{1/3})}$. This gives the statement of the lemma.\sp

\noi{\bf Lemma 2}. We have as $m\ra\iy$
\[m^{1/3}\,\ps(u)\inv\,v^{-m^{1/3}x}\,{1\ov 2\pi i}\int
\ps(z)\,z^{m^{1/3}\,x}\,dz=
-(3b)^{-1/3}{\rm Ai}(x/v(3b)^{1/3})+O(m^{-1/3})\]
uniformly for bounded $x$.\sp

\noi{\bf Proof}. Near $z=u$
\[\si(z)=\si(u)+ b(z-u)^3+O((z-u)^4)),\]
\[\log (-z)=\log v+{1\ov u}(z-u)+O((z-u)^2)=\log v-{1\ov
v}(z-u)+O((z-u)^2).\]
Hence, using Lemma 1, we have
\[{1\ov 2\pi i}\int_{C^+}\ps(z)\,z^{m^{1/3}\,x}\,dz\]
\[=O(|\ps(u)|e^{-\dl m})+\ps(u)\,v^{m^{1/3}x}
{1\ov 2\pi i}\int_{C^+_\ep}e^{m b(z-u)^3-m^{1/3}x(z-u)/v
+O(m(z-u)^4+m^{1/3}x(z-u)^2)}\,dz.\]

We show that removing the $O$ term in the exponential in the integrand
leads to an error $O(m^{-2/3})$ in the integral. This error equals
\[\int_{C^+_\ep}e^{mb(z-u)^3-m^{1/3}x(z-u)/v}\left(
e^{O(m(z-u)^4+m^{1/3}x(z-u)^2)}-1\right)\,dz\]
\[=\int_{C^+_\ep}e^{mb(z-u)^3-m^{1/3}x(z-u)/v+
O(m(z-u)^4+m^{1/3}x(z-u)^2)}\,
O(m(z-u)^4+m^{1/3}x(z-u)^2)dz.\]
Now the exponential has the form
\[e^{m b(1+\eta_1)(z-u)^3-m^{1/3}x(1+\eta_2)(z-u)/v},\]
where the $\eta_i$ can be made arbitrarily small by taking $\ep$ small
enough.
If we make the substitution $z-u=\z/m^{1/3}$ the error becomes
\[m^{-2/3}\int e^{b(1+\eta_1)\z^3-x(1+\eta_2)\z/v}\,O(\z^4+x\z^2)\,d\z.\]
The integral is now taken over a long contour lying in thin angles
around the rays $|{\rm arg}\,\z|=\pi/3$, with ends having absolute value
at least a constant times $m$. This integral is clearly bounded,
uniformly
in $m$ for bounded $x$.

Therefore with the stated error we may remove the $O$ terms from the
exponential in the original integral. Then we make the same
substitution.
The integrand is exponentially small at the ends of the resulting
contour. Therefore if we complete it so that it goes to infinity in
the two directions $\pm \pi/3$ the error incurred will be exponentially
small.

We have shown that
\[m^{1/3}\,\ps(u)\inv\,v^{-m^{1/3}x}\,{1\ov 2\pi
i}\int\ps(z)\,z^{m^{1/3}\,x}\,dz=
{1\ov 2\pi i}\int_{\iy e^{i\pi/3}}^{\iy e^{-i\pi/3}}
e^{b\z^3-x\z/v}d\z+O(m^{-1/3}).\]
If we deform the contour to the imaginary axis and make the substitution
$\z\ra -i\z$ then the last integral, with
its factor, becomes
\[-{1\ov 2\pi }\int_{-\iy}^{\iy }e^{ib\z^3+ix\z/v}d\z=-(3b)^{-1/3}{\rm
Ai}(x/v(3b)^{1/3}).\]
This proves the lemma.\sp

The second integral in (\ref{psints}) is similar.\sp

\noi{\bf Lemma 3}. We have as $m\ra\iy$
\[m^{1/3}\,\ps(u)\,v^{m^{1/3}x}\,{1\ov 2\pi i}\int
_{C^-}\ps(z)\inv\,z^{-m^{1/3}\,x}\,dz=
-(3b)^{-1/3}{\rm Ai}(x/v(3b)^{1/3})+O(m^{-1/3})\]
uniformly for bounded $x$.\sp

\noi{\bf Proof}. The derivation is essentially the same. The
exponentials are
replaced by their negatives and the directions $\pm \pi/3$ are replaced
by $\pm 2\pi/3$.
The fact that $C^-$ is unbounded causes no difficulty since the
integrand now behaves at
infinity like a large negative power of $z$.
We get the same Airy function in the end, as we have already seen.\sp

Now for the tricky part. We need estimates that are uniform for all $x$
and where the
error term contains a factor which is very small for large $x$. In fact
we shall show
that the statements of Lemmas 2 and 3 hold, uniformly for all $x$, when
the error terms
are replaced by $m^{-1/3}e^{-x}$. (The $x$ in the exponential can be
improved to
a constant times $x^{3/2}$ but that makes no difference.) To do this we
have to be
more careful and use the steepest descent curves for the full integrands
in
(\ref{psints}),
not just for the factors $\ps^{\pm 1}$. We consider in detail only the
first
integral in (\ref{psints}); as before, the second is treated
analogously.

Set
\[\ps(z,c')=(1+z)^n\,(r-z)^m\,(-z)^{-(1-c')\,m}=\ps(z)\,(-z)^{(c'-c)\,m}.\]
We are interested in the asymptotics of
\be I(c')={1\ov2\pi i}\int \ps(z,c')\,dz\label{Iint}\ee
when $c'-c=m^{-2/3}\,x$. Our condition on $x$ says that $c'<1$ and, in
view of
what we already know, we may assume $x$ is positive and bounded away
from zero, so $c'>c$.

We let $C$ be the steepest descent curve for
$\ps(z,c')$. This curve now
passes vertically
through one of the critical points of $\ps(z,c')$. For $c'>c$ there are
two critical
points
\[u^{\pm}_{c'}={-(1-r)c'-(1-\al)r \pm
\sqrt{\left((1+r)\,c'+(\al-1)\,r\right)^2-4\al\, r}\ov2(\al+c')},\]
which are real and satisfy
\[-1<u^-_{c'}<-v<u^+_{c'}<0.\]
To determine which critical point our curve passes through we consider
the function
\[\si(z,c')={1\ov m}\log\,\ps(z,c')=\si(z)+(c'-c)\log\,(-z).\]
The critical points $u^{\pm}_{c'}$
are the zeros of $\si_z(\u,c')$. (Subscripts here and below denote
derivatives in the
usual way.)
We use the fact that $u^{\pm}_{c'}$ are smooth functions of
$\ga=\sqrt{c'-c}$ and compute, recalling that
$\si_z(u,c)=\si_{zz}(u,c)=0$ and observing that $dc'/d\ga=0$ when
$\ga=0$,
\be{d\ov d\ga}\si_{zz}(u^{\pm}_{c'},c')\Big|_{\ga=0}=\si_{zzz}(u,c)
{du^{\pm}_{c'}\ov d\ga}\Big|_{\ga=0}.\label{sizz}\ee
The first factor on the right is positive (we denoted it by $6b$), while
\be{du^{\pm}_{c'}\ov d\ga}\Big|_{\ga=0}=\pm\beta\label{dudga}\ee
where
\[\beta={(\al r)^{1/4}\,(1+r)^{3/2}\ov (\sqrt{\al}+\sqrt{r})^2} .\]
Since $\si_{zz}(u^{\pm}_{c'},c')=0$
when $\ga=0$ we deduce that for $c'$ close to, but greater than, $c$ we
have
\[\si_{zz}(u^+_{c'},c')>0,\ \ \ \si_{zz}(u^-_{c'},c')<0.\]
These inequalities hold for all $c'$ since the second
derivative can be zero only when $c'=c$. This shows that the steepest
descent curve $C$
for $\ps(z,c')$ passes through $u^+_{c'}$, because on the curve
$|\ps(z,c')|$ has a
maximum at the critical point. (Similarly the steepest descent curve for
$\ps(z,c')\inv$ passes through $u^-_{c'}$.) To make the notation less
awkward we write
$u_{c'}$ instead of $u^+_{c'}$. First, we have the analogues of Lemmas 2
and 3.\sp

\noi{\bf Lemma 4}. Given $\ep>0$ there exists a $\dl>0$ such that
$I(c')=O(|\ps(u,c')|e^{-\dl\,m})$ if $c'-c>\ep$.\sp

\noi{\bf Proof}. The function $\si(z,c')$ is decreasing for $u<z<\u$
since it decreases
near and to the left of $\u$ and has no critical point in this interval.
Hence
$\si(\u,c')<\si(u,c')$ so
$\si(\u,c')-\si(u,c')$ is negative and bounded away from zero for
$c'>c+\ep$.
Since the maximum of $|\ps(z,c')|$ on $C$ is at $z=\u$ the statement
follows.\sp

In view of Lemma 4 we may assume in what follows that $c'-c$ is as small
as we please.
We denote by $C_\ep$ the portion of $C$ which lies within $\ep$ of the
critical point $\u$.\sp

\noi{\bf Lemma 5}. If in the integral (\ref{Iint}), in which we
integrate over
$C$, we integrate only over $C_\ep$ the error incurred
is $O(|\ps(u,c')|\,e^{-\dl\,m})$ for some $\dl>0$. \sp

\noi{\bf Proof}. The maximum of $|\ps(z,c')|^{1/m}$ on $C$ occurs at
$\u$ and it is stricly
smaller on the complement of $C_\ep$ in $C$ than it is at $\u$.
Therefore
the integral in question is $O(|\ps(\u,c')|\,e^{-\dl\,m})$ for some
$\dl>0$.
Since $\si(\u,c')<\si(u,c')$, as we saw in the proof of the last lemma,
this one
is established.\sp

Because of Lemmas 4 and 5 we need only compute the behavior of
$\ps(z,c')$, or
equivalently $\si(z,c')$, for $z$ near $\u$. Recall that $u=u_{c}=-v$.
\sp

\noi{\bf Lemma 6}. We have
\begin{eqnarray*}
(i)\ \si(\u,c')&=&\si(u)+(c'-c)\,\log v-{2\ov3}{\beta \ov
v}\,(c'-c)^{3/2}+O\left((c'-c)^2\right)\;;\\
(ii)\ \si_z(\u,c')&=&0\;;\\
(iii)\ \si_{zz}(\u,c')&=& 6b\beta\,\sqrt{c'-c}+O(c'-c)\;;\\
(iv)\ \si_{zzz}(\u,c')&=& 6b+O(\sqrt{c'-c}\,)\;.
\end{eqnarray*}\sp

\noi{\bf Proof}. From (\ref{dudga}) and the fact that $\u$ is a smooth
function of $\ga$
(or directly) we see that
\be\u-u=\beta\sqrt{c'-c}+O(c'-c).\label{u1}\ee
Consequently, since $u=-v$,
\be {\u\ov u}=1-{\beta\ov v}\sqrt{c'-c}++O(c'-c).\label{u2}\ee
Now since  $\si_z(\u,c')=0$ we have
\[{d\ov dc'}\si(\u,c')=
{\partial\ov \partial c'}\si(z,c')\Big|_{z=\u}=\log(-\u)=\log
v+\log{\u\ov
u}.\]
Integrating with respect to $c'$ from $c$ to $c'$ and using (\ref{u2})
we obtain (i).
Of course (ii) is immediate. As for (iii) and (iv), these follow from
(\ref{sizz})
and (\ref{dudga}) and the fact that $\u$ is a smooth functions of
$\ga$.\sp

\noi{\bf Lemma 7}. The conclusions of Lemmas 2 and 3 hold uniformly for
all $x$
when the error terms are replaced by $O(e^{-\dl m})+O(m^{-1/3}e^{-x})$
for some $\dl>0$.\sp

\noi{\bf Proof}. We consider (\ref{psint}), which is $I(c')$ with
$c'-c=m^{-2/3}x$.
Putting together Lemmas 5 and 6 we deduce that
\[I(c')=O(|\ps(u)\,v^{(c'-c)m}|\,e^{-\dl\,m})+\ps(u)\,v^{(c'-c)m}\,
e^{-{2\beta\ov 3v}(c'-c)^{3/2}m}\times\]
\[{1\ov 2\pi
i}\int_{C_\ep}e^{mb(z-\u)^3+3mb\beta\sqrt{c'-c}\,(z-\u)^2+O(m[(c'-c)^2+
(c'-c)|z-\u|^2+\sqrt{c'-c}|z-\u|^3+|z-\u|^4])}\,dz.\]
If $c'-c=m^{-2/3}x$ the exponential factor equals $e^{-{2\beta\ov
3v}x^{3/2}}$
while the integral equals
\[{1\ov 2\pi i}\int_{C_\ep}e^{m b(z-\u)^3+3m^{2/3}bx^{1/2}\beta (z-\u)^2
+O(m^{-1/3}x^2+m^{1/3}x|z-\u|^2+m^{2/3}x^{1/2}|z-\u|^3+m|z-\u|^4])}\,dz.\]
Now $C_{\ep}$, rather than looking like two rays near the critical
point, looks
like one branch of a hyperbola.

Note that by Lemma 4 we may assume that
$c'-c=m^{-2/3}x$ is as small as desired. It follows that the
exponent, without the $O(m^{-1/3}x^2)$ term, can be written
\[m(b+\eta_1)\,(z-\u)^3+3m^{2/3}(b+\eta_2)x^{1/2}\beta\,(z-\u)^2,\]
where, if $\ep$ is chosen small enough, the $\eta_i$ can be made as
small as
desired. Upon making the variable change $z-\u=\z/m^{1/3}$
the integral becomes
\[{m^{-1/3}\ov 2\pi i}\int
e^{(b+\eta_1)\z^3+3(b+\eta_2)x^{1/2}\beta\z^2}\,d\z,\]
taken over a long contour in the right half-plane on which
$|{\rm arg}\,\z|>\pi/3-\eta$,
with another small $\eta$. The integral here is uniformly bounded.

To take care of the term $O(m^{-1/3}x^2)$ in the exponential in the
original integral,
 observe that
if $m^{-2/3}x$ is small enough then
$m^{-1/3}x^2$ will be at most a small constant times $x^{3/2}$ and
so
\[e^{-{2\beta\ov 3v}x^{3/2}}\left(e^{O(m^{-1/3}x^2)}-1\right)
=O(m^{-1/3}x^2\,e^{-{\beta\ov 2v}x^{3/2}})=O(m^{-1/3}e^{-x}).\]
Thus removing the term from the exponential leads to an eventual error
$O(m^{-2/3}e^{-x})$. That removing the other $O$ terms from the
exponential
leads to the same error is seen as it was in the proof of Lemma 2---the
substitution
in the integral representating the error results in an extra factor
$m^{-1/3}$ and there is the exponential factor $e^{-{2\beta\ov
3v}x^{3/2}}$
outside the integral.

After removing all the $O$ terms and making the variable change
$z-\u=\z/m^{1/3}$
the integral becomes
\[{m^{-1/3}\ov 2\pi i}\int e^{b\z^3+3bx^{1/2}\beta\z^2}\,d\z,\]
taken over a long contour in the right half-plane on which
$|{\rm arg}\,\z|>\pi/3-\eta$.  Completing the
contour so that it goes to infinity in the directions ${\rm arg}\,\z=\pm
\pi/3$ leads to an exponentially small error. It follows that (the first
part
of) the lemma
holds with the negative of the Airy function in the statement replaced
by
\[e^{-{2\beta\ov 3v}x^{3/2}}{1\ov2\pi i}
\int_{\iy e^{i\pi/3}}^{\iy
e^{-i\pi/3}}e^{b\z^3+3bx^{1/2}\beta\z^2}d\z.\]
If we complete the cube and make the substitution $\z\ra\z-\beta
x^{1/2}$
this becomes, upon noting that $3 b \beta^2 =1/v$,
\[\int_{\iy e^{i\pi/3}}^{\iy
e^{-i\pi/3}}
e^{b\z^3-x\z/v}d\z
=-(3b)^{-1/3}{\rm Ai}\,(x/(v(3b)^{1/3})).\]

The second part of the lemma is analogous, just as the proof of Lemma 3
was analogous to the proof of Lemma 2.\sp

We have now shown that if we set $j=m^{1/3}x,\ k=m^{1/3}y$ then
\[(-1)^{h+j+k}m^{1/3}\,\ps(u)\,v^{m^{1/3}s+j+k}(\ph_+/\ph_-)_{-h-j-k-1}\ra
-(3b)^{-1/3}{\rm Ai}((s+x+y)/v(3b)^{1/3}),\]
and the difference between the two is $O(m^{-1/3}e^{-(x+y)})+O(e^{-\dl
m})$.
It follows easily from this that if we denote the matrix on the left,
without the
factor $m^{1/3}$, by $(M(j,k))$ and the
kernel on the right by $A(x,y)$ then the kernel
$m^{1/3}\,M([m^{1/3}\,x],\,[m^{1/3}\,y])$
converges in Hilbert-Schmidt norm to the kernel $A(x,y)$ on $(0,\,\iy)$.
(Recall that $j$ and $k$
are at most $O(m)$. Therefore the error term $O(e^{-\dl m})$ can only
contribute an
exponentially small error to the norm and so can be ignored.
Similarly
we can let our indices $j$ and $k$ run to infinity.)
Thus, under the scaling $j\ra m^{1/3} x,\ k\ra m^{1/3} y$ the
matrices with $j,k$ entry
\[(-1)^{h+j+k}\,\ps(u)\,v^{m^{1/3}s+j+k}(\ph_+/\ph_-)_{-h-j-k-1}\]
scale in Hilbert-Schmidt norm to the kernel
$A(x,y)$. Similarly
so do the matrices with $j,k$ entry
\[(-1)^{h+j+k}\,\ps(u)\inv\,v^{-(n^{1/3}s+j+k)}(\ph_-/\ph_+)_{h+j+k+1}.\]
Therefore the product of the matrices scale in trace norm to the
(operator)
square of the kernel, which is the Airy kernel
(\ref{KAiry}). This or completeness
 the justification.

\subsection{Critical Regime: $p\sim p_c$}
When $p=p_c$ ($\al r=1$),\footnote{See the remark at
end of this section.}
 the analysis of the previous section must be
modified.
  We set $h=m-\h$ ($\h=0,1,2,\ldots$) and introduce the new
$\ps$
\[ \ps=(1+z)^n\, (z-r)^m \]
and the corresponding new $\s$
\[ \s(z)={1\ov m}\,\log\ps = \al\log(1+z)+\log(z-r). \]
The saddle point now occurs at $z=0$ with
$\s^{\prime\prime}(0)=-\al(1+\al)$.
Thus in the neighborhood of $z=0$
\be \ps(z) \sim (-1)^m r^m \,e^{-m\al(1+\al)\, z^2/2}. \label{psi0}\ee

Since $(\vp_+/\vp_-)_{-h-k-j-1}$ vanishes for $h+j+k+1>m$, we can again
assume
$h+j+k<m$ which becomes the condition $j+k<\h$.  As before our starting
point is
 the integral expression
\[ (\vp_+/\vp_-)_{-h-j-k-1}={1\ov 2\pi i}\int \ps(z) z^{-m+h+j+k}\, dz
\]
where the contour is a circle centered at 0 with radius $\rho<1$.
Taking
this $\rho$ sufficiently
small so that we may use the approximation (\ref{psi0}) on the
integrand,
we obtain  after making
the change of variables
\[ \z=\left({m\al(1+\al)\ov 2}\right)^{1/2} \, z = z/S,\]

\ba (\vp_+/\vp_-)_{-h-j-k-1} &\sim & (-1)^m\, r^m \,S^{j+k-\h+1}
\,{1\ov 2\pi i}\, \int e^{-\z^2} \z^{j+k-\h}\, d\z \\
&=& \left\{\begin{array}{cc}
(-1)^m\, r^m \,S^{j+k-\h+1}\> {(-1)^L\ov L!} & \textrm{if}\>\>
\h-j-k-1=2L=0,2,4\ldots \\
0 & \textrm{if} \>\> \h-j-k-1=\textrm{odd integer}.
\end{array}\right.
\ea

Our second integral is
\[ (\vp_-/\vp_+)_{h+j+k+1}={1\ov 2\pi i}\, \int\ps(z)\inv
z^{m-h-j-k-2}\, dz \]
where the contour has $-1$ on the outside and $r$ on the inside.
We deform the contour to the imaginary axis going from $i\iy$ to
$-i\iy$ with an infinitesimal indentation going around
0 to the left.  The part of the contour lying in the right
half plane is exponentially small because of the factor $(1+z)^{-n}$ and
can therefore be neglected.
For the  integral along
the imaginary axis we can replace $\ps$ by (\ref{psi0}) with
an error that is exponentially small.  Thus the above integral
is asymptotically equal  to
\[(-1)^m {r^{-m}\ov 2\pi i}\, \int_{i\iy}^{-i\iy} e^{z^2/S^2}
z^{\h-j-k-2}\, dz, \]
which in turn equals
\[(-1)^m i^{\h-j-k-1}r^{-m}\, S^{\h-j-k-1}\,
{1\ov 2\pi i}\int_{\iy}^{-\iy} e^{-\z^2} \z^{\h-j-k-2}
\, d\z \]
where there is an indentation above $\z=0$.
If we now substitute $\z=\sqrt{t}$, the above integral becomes
\[ (-1)^mi^{\h-j-k-1}r^{-m}\, S^{\h-j-k-1}\, {1\ov 4\pi i}
\int_{\iy}^{0^{+}} e^{-t} t^{(\h-j-k-1)/2-1}\,dt. \]
The contour starts at $+\iy$, loops around 0 in the positive
direction and then returns to $+\iy$.  This last integral is
Hankel's integral representation of the $\Gamma$ function. Thus
\[ (\vp_-/\vp_+)_{h+j+k+1}\sim {(-1)^{h+j+k+1}\ov 2\pi}\,r^{-m}\,
 S^{\h-j-k-1}\,\sin\left({\pi\ov 2}(\h-j-k-1)\right)\,
 \Gamma\left({\h-j-k-1\ov 2}\right).\]

 We now use these two asymptotic expressions along with the condition
 $j+k<\h$ in (\ref{Kdefn}) to obtain (after a short calculation)
 \be K_h(j,k)\sim {(-S)^{k-j}\ov 2\pi}\, \sum_{\ell=0}^{[{\h-k-1\ov 2}]}
 {1\ov \ell!}\,\sin{\pi\ov 2}(k-j)\,\Gamma\left(\ell+{k-j\ov 2}\right).
\label{Kcrit}\ee
When $\ell+(k-j)/2$ is a nonpositive integer, the product of the sine
and gamma functions
is replaced by
\[ {(-1)^\ell\, \pi\ov \left({j-k\ov 2}-\ell\right)!}\, . \]
The factor $(-S)^{k-j}$ may be dropped when computing the determinant
$\det(I-K_h)$ since
it does not change its value.  We evaluate this determinant and display
the results
for $\h\le 9$ in Table~\ref{critProbTable}.

\textit{Remark}. Since $m$ and $n$ are integers it is extremely
unlikely that $p=p_c=m/(m+n)$.  If $p$ is irrational this never
occurs.  However the preceding analysis shows that if $\al r=1+o(m^{-1/2})$
rather than 1 then in the integrals one gets extra factors
$(1+z)^{o(m^{1/2})}$.  Then after the substitution $z=S\zeta$ this drops
out since $S=O(m^{-1/2})$.  The upshot is that the asymptotics hold for
any $p$ when $m$ and $n$ go to infinity in such a way that
$m/(m+n)=p+o(m^{-1/2})$.

\begin{table}
\bc
\small{
\begin{tabular}{|r|l|l|}\hline
$\h$ & $\lim_{m\ra\iy}\pr\left(H-m\le -\h\right)$ & Numerical Value \\
\hline
0 & $1$ & 1.0  \\[1ex]
1 & ${1\ov 2}$ & 0.5  \\[1ex]
2 & ${1\ov 4}-{1\ov 2\pi}$  & $9.08451\times 10^{-2}$ \\[1ex]
3 & ${1\ov 8} -{3\ov 8\pi}$  & $5.63379\times 10^{-3}$\\[1ex]
4 & ${1\ov 16}+{1\ov 3\pi^2}-{29\ov 96\pi}$ & $1.17616\times 10^{-4}$
\\[1ex]
5 & ${1\ov 32}+{41\ov 144\pi^2}-{145\ov 768\pi}$  & $8.22908 \times
10^{-7}$ \\[1ex]
6 & ${1\ov 64}-{32\ov 135\pi^3}+{1169\ov 3840\pi^2}-{1249\ov 10240\pi}$
& $1.92570\times 10^{-9}$ \\[1ex]
7 & ${1\ov 128}-{49\ov 225\pi^3}+{198827\ov 921600\pi^2}-{8743\ov
122880\pi}$ & $1.50565\times 10^{-12}$\\[1ex]
8 & ${1\ov 256}+{4096\ov 23625\pi^4}-{10289\ov 36000\pi^3}+{5773487\ov
34406400\pi^2}-{145603\ov 3440640\pi}$
 & $3.92048\times 10^{-16}$\\[1ex]
 9 & ${1\ov 512} +{15376\ov 91875\pi^4}-{5528469\ov
25088000\pi^3}+{279234531\ov 2569011200\pi^2}-
 {436809\ov 18350080\pi}$& $3.42524\times 10^{-20}$\\[1ex]
\hline
\end{tabular}
}
\caption{Limiting Distribution when $p\sim p_c$ \label{critProbTable}}
\ec
\end{table}

\subsection{Deterministic Regime: $p>p_c$}
\subsubsection{Large Deviations Approach}
Assume that $p>p_c$. Then there exists an $\e>0$ so
that $n/m$ approaches $(1+\e)(1/p-1)$. To simplify
the statements, we will just assume that $n=(1+\e)(1/p-1)m$.

Imagine the random $m\times n$ matrix $A$ from \S2.1 as the lower
left corner of an infinite matrix of 0's and 1's, created
by the independent coin flips.
Fix a position $(i,j)$ ($i, j\ge 1$)
in this infinite random matrix.
Define $J$ as the column index
of the first entry, from left to right, with a 1 on the row
{\it above\/} $(i,j)$ and in the columns larger or equal $j$.
Then define $\xi_{(i,j)}=J-j$.
In the example given, $\xi_{(3,1)}=0$ and $\xi_{(5,1)}=3$.

Now create a sequence of i.i.d.\  random variables
$\xi_1, \xi_2,\dots$, as follows.
Let $\xi_1$ equal the
column index minus one of the first 1 on the first row.
Then let $\xi_2=\xi_{(1,1+\xi_1)}$, $\xi_3=\xi_{(2,1+\xi_1+\xi_2)},
\dots$.
The basic observation is that, since we are always taking the
best positioned 1 on the next line, we have equality of the two events
\[
\{\textrm{there is an increasing path of length $m$ in $A$}\}=
\{\xi_1+\dots+\xi_m<n\}.
\]
Therefore, we need to show that
\[
\pr(\xi_1+\dots+\xi_m\ge n)\]
goes to 0 exponentially as $m\to\infty$. However,
$\pr(\xi_1=i)=p(1-p)^i$,
$i=0,1,\dots$ and so $E(\xi_1)=1/p-1$. By elementary
large deviations (e.g.\ \S1.9 in \cite{durrett}),
\[
-m^{-1}\log P(\xi_1+\dots+\xi_m\ge n)\to \ga(\e)\]
where an elementary calculation shows 
\[ \ga(\e)=(1/p-1)(1+\e)\log(1+\e)
-p^{-1}(1+\e-\e p)\log(1+\e-\e p)\]
which is positive whenever $\e>0$.

\subsubsection{Saddle Point Approach}
For completeness, we show how the saddle point
method gives the same result.
Thus we show that when $p>p_c$ (or $\al r>1$)
\[\det(I-K_h)\ra0\]
exponentially as $m\ra\iy$ even when $h=m-1$, thus establishing
assertion
1(c) in \S3  with exponential approach to the limit.

As we saw at the beginning in the last section we need only consider the
entries
$K_h(j,k)$ when $h+j+k<m$, which in the present situation means $j=k=0$.
Our
claim is therefore that $K_{m-1}(0,0)\ra 1$ exponentially as $m\ra\iy$.
The first
integral to consider is
\[ (\vp_+/\vp_-)_{-m}={1\ov 2\pi i}\int \ps(z)\, z^{-1}\, dz
=\ps(0)=(-r)^m.\]
 The second integral is
\[ (\vp_-/\vp_+)_{m}={1\ov 2\pi i}\, \int\ps(z)\inv z^{-1}\, dz .\]
Recall that the contour here surrounds 0 and has $-1$ on the outside,
$r$ on the inside.
The critical point for steepest descent is at $z=u$ where
\[{\al\ov 1+u}+{1\ov u-r}=0,\ \ \ \ \ u={\al r-1\ov \al+1}.\]
The steepest descent curve will pass vertically through this point and
go to $\iy$
in two directions. But notice that since $u$ is positive, in order to
deform our original
contour to this one we have to pass through $z=0$. The residue of the
integrand there
equals $(-r)^{-m}$ and so
\[ (\vp_-/\vp_+)_{m}=(-r)^{-m}+{1\ov 2\pi i}\int \ps(z)\inv z^{-1}\,
dz,\]
where now the integral is taken over the steepest descent curve. This
integral
is asymptotically a constant times $m^{-1/2}$ times the value of the
integrand at $z=u$,
and this value equals $(-1)^m$ times
\[\left({\al\,(r+1)\ov \al+1}\right)^{-\al m}\left({r+1\ov
\al+1}\right)^{-m}.\]
Our claim is therefore equivalent to the statement that this is
exponentially smaller than
$r^{-m}$, which in turn is equivalent to the inequality
\[(r+1)^{\al+1}{\al^{\al}\ov(\al+1)^{\al+1}}>r.\]
It is an elementary exercise that this is true for all $r\ge0$ except
for $r=1/\al$, when
equality holds. But in our case $r>1/\al$ so the inequality holds.

\subsection{Finite GUE Regime: Fixed $x$ and $t\ra\iy$}
\subsubsection{Saddle Point Calculation}
We return to (\ref{phi+-}) and this time set
\[h={r\ov 1+r}m+s\,m^{1/2}=p\,m+s\,m^{1/2},\]
and make the substitutions
$j\ra x\,m^{1/2},\;k\ra y\,m^{1/2}$
to write the integral (\ref{phi+-}) as
\be{1\ov 2\pi
i}\int(1+z)^n\,(r-z)^m\,(-z)^{-m/(1+r)}\,(-z)^{(s+x+y)\,m^{1/2}}\,dz.
\label{int1}\ee
Now we set
\[\ps(z)=(r-z)^m\,(-z)^{-m/(1+r)},\]
which is the main part of the integrand. There is a single critical
point, $z=-1$, and
at this point $d^2/dz^2\,\log\ps(z)$ is equal to
\[m{r\ov(1+r)^2}=m\,p\,(1-p).\]
This is positive and so the steepest descent curve is vertical at the
critical
point; it goes around the origin and closes at $z=r$. The main
contribution
to the integral comes from the immediate  neighborhood of the critical
point. If we make
the variable change
\[z=-1+{\z\ov\sqrt m}\]
and take into account the other factors in the integrand we see the
integral is
asymptotically
\be -{(r+1)^m\ov2\pi i}\int_{-i\iy}^{i\iy}\left({\z\ov\sqrt
m}\right)^n\,
e^{{1\ov2}{p(1-p)}\z^2-(s+x+y)\z}\,{d\z\ov\sqrt
m}.\label{int2}\ee

Now
\[{1\ov2\pi
i}\int_{-i\iy}^{i\iy}e^{a\z^2-b\z}\,d\z={e^{-b^2/4a}\ov2\sqrt{a\pi}}\]
and so
\[{1\ov2\pi
i}\int_{-i\iy}^{i\iy}\z^n\,e^{a\z^2-b\z}\,d\z={(-1)^n\ov2\sqrt{a\pi}}
{d^n\ov db^n}e^{-b^2/4a}={1\ov\sqrt\pi(2{\sqrt a})^{n+1}}\,
e^{-b^2/4a}\,H_n\left({b\ov2\sqrt a}\right).\]
($H_n$ are the Hermite polynomials.)
Hence our first integral (\ref{int1}) is asymptotically equal to
 $-(r+1)^m/\sqrt
m^{n+1}$ times this
expression with
\[a={1\ov2}\,p\,(1-p),\ \ \ \ b=s+x+y.\]
Thus we have shown that the matrix with $j,\,k$ entry
\[(-1)^{-h-j-k}(r+1)^{-m}\,m^{n/2}\,(\ph_+/\ph_-)_{-h-j-k-1}\]
scales to the operator on $(0,\,\iy)$ with kernel
\[-{1\ov\sqrt\pi(2{\sqrt
a})^{n+1}}\,e^{-(s+x+y)^2/4a}\,H_n\left({s+x+y\ov2\sqrt a}\right),\]
with $a$ as given above.

Next, with the same substitutions in the integral,
\[(\ph_-/\ph_+)_{h+j+k+1}={1\ov 2\pi
i}\int(1+z)^{-n}\,(z-r)^{-m}\,z^{m-h-j-k-2}\,dz\]
\[=(-1)^{-h-j-k}{1\ov 2\pi i}\int(1+z)^{-n}\,(r-z)^{-m}\,
(-z)^{m/(1+r)}\,(-z)^{-(s+x+y)\,m^{1/2}-2}\,dz.\]
The contour here encloses $0$ and $r$ and has $-1$ on the outside. The
steepest
descent curve should go through the critical point $-1$ horizontally. We
deform the
given contour to a curve starting at $-\iy+0i$,
going above the the real axis, looping around $z=-1$ clockwise, then
back
below the real axis
to $-\iy-0i$. The original contour can be deformed to this because the
integrand is
small at $\iy$. The main
contribution is again in the neighborhood of $z=-1$. Making the same
variable change
as before leads to an integral which is asymptotically
\be -{(r+1)^{-m}\ov2\pi i}\int\left({\z\ov\sqrt m}\right)^{-n}\,
e^{-{1\ov2}{p(1-p)}\z^2+(s+x+y)\z}\,{d\z\ov\sqrt
m},\label{int3}\ee
where now the contour is a circle going around $\z=0$ counterclockwise.
Using now the fact
\[{1\ov2\pi
i}\int\z^{-n}\,e^{-a\z^2+b\z}\,d\z={a^{(n-1)/2}\ov(n-1)!}e^{-b^2/4a}\,
{d^{n-1}\ov d\z^{n-1}}e^{-(\z-b/2\sqrt a)^2}\Big|_{\z=0}\]
\[={a^{(n-1)/2}\ov(n-1)!}\,H_{n-1}\left({b\ov2\sqrt a}\right)\]
we find that the matrix with $j,\,k$ entry
\[(-1)^{h+j+k}(r+1)^{m}\,m^{-n/2}\,(\ph_-/\ph_+)_{h+j+k+1}\]
scales to the operator on $(0,\,\iy)$ with kernel
\[-{a^{(n-1)/2}\ov(n-1)!}\,H_{n-1}\left({s+x+y\ov2\sqrt a}\right).\]

Combining, we see that the product of the two matrices (aside from a
factor
$(-1)^{j-k}$, which does not affect the determinant) has scaling limit
the operator with kernel
\[{1\ov\sqrt\pi2^{n+1}a\,(n-1)!}\int_0^{\iy}e^{-(s+x+z)^2/4a}\,
H_n\left({s+x+z\ov2\sqrt a}\right)\,H_{n-1}\left({s+z+y\ov2\sqrt
a}\right)\,dz.\]

Instead of a direct evaluation of this last integral, we
will not evaluate our $\z$ integrals (\ref{int2}) and (\ref{int3}), but
rather consider
them as integrals with variables
$\z_1$ and $\z_2$, combine and integrate with respect to $z$. We see
that the scaled
kernel for the product is
\[-{1\ov 4\pi^2}\int_0^{\iy}\int\int\left({\z_1\ov\z_2}\right)^n
e^{a(\z_1^2-\z_2^2)-(s+x+z)\z_1
+(s+z+y)\z_2}\,dz\,d\z_1\,d\z_2,\]
where the $\z_1$ contour is a vertical line described upward and the
$\z_2$
contour goes around 0 counterclockwise. If the vertical line is to the
right
of the circle we can integrate first with respect to $z$, yielding
\[-{1\ov
4\pi^2}\int\int\left({\z_1\ov\z_2}\right)^ne^{a(\z_1^2-\z_2^2)-(s+x)\z_1
+(s+y)\z_2}\,{d\z_1\,d\z_2\ov \z_1-\z_2}.\]

Let's call this $L_n(x,\,y)$. This is 0 when $n=0$, and
\[L_k(x,\,y)-L_{k-1}(x,\,y)=-{1\ov 4\pi^2}\int\int{\z_1^{k-1}\ov\z_2^k}
e^{a(\z_1^2-\z_2^2)-(s+x)\z_1
+(s+y)\z_2}\,d\z_1\,d\z_2.\]
This integral is a product and we can use the computations we did above
to see that
it equals
\[{1\ov 2^k\sqrt{\pi
a}(k-1)!}e^{-(s+x)^2/4a}\,H_{k-1}\left({s+x\ov2\sqrt a}\right)\,
H_{k-1}\left({s+y\ov2\sqrt a}\right).\]
If $\ph_k$ are the oscillator wave
functions\footnote{The
oscillator wave functions are $\vp_k(x):=e^{-x^2/2}\,H_k(x)/
\sqrt{2^k k!\, \pi^{1/2}}$
and form an orthonormal basis for $L^2((0,\iy))$.
The Hermite kernel is $K_{H,n}(x,y):=\sum_{k=0}^{n-1}\ph_k(x)\ph_k(y)$.}
then this equals
\[{1\ov2\sqrt a}\,\ph_{k-1}\left({s+x\ov2\sqrt a}\right)\,
\ph_{k-1}\left({s+y\ov2\sqrt a}\right)\]
times the factor
\[e^{-(s+x)^2/8a}\,e^{(s+y)^2/8a}.\]
It follows that if $K_{H,n}$ is the Hermite kernel then
\be L_n(x,\,y)=e^{-(s+x)^2/8a}\,{1\ov2\sqrt
a}\,K_{H,n}\left({s+x\ov2\sqrt
a},
\,{s+y\ov2\sqrt a}\right)\,e^{(s+y)^2/8a}.\label{LK}\ee
We deduce that
\[ \lim_{m\ra\iy}{\rm Prob}\;(H\le p\,m+ s\,m^{1/2})\]
is equal to the Fredholm determinant of
\[{1\ov2\sqrt a}\,K_{H,n}\left({s+x\ov2\sqrt a},\,{s+y\ov2\sqrt
a}\right)\]
over $(0,\,\iy)$, or equivalently the Fredholm determinant of
$K_{H,n}(x,\,y)$ over $(s/2\sqrt a,\,\iy)$.  It is notationally
convenient to introduce
\[ \s^2:=2 a = p(1-p) \]
and to define
\[ F_n^{GUE}(s):=\lim_{m\ra\iy}\pr\left({H-p m\ov \s \sqrt{m}}\le
s\right). \]
This equals  the Fredholm determinant of $K_{H,n}$ over
$(s/\sqrt{2},\iy)$
and is equal to the distribution of the largest eigenvalue in the finite
$n$ GUE.\footnote{Our
normalization of $F_n^{GUE}$ differs from the usual one \cite{mehta,tw3}
by
a factor of $\sqrt{2}$, i.e.\ the usual normalization is the Fredholm
determinant of the Hermite
kernel over $(s,\iy)$.}

\subsubsection{Moments of $F_n^{GUE}$}

{}From the theory of random matrices, e.g.~\cite{mehta,tw3},
 we know that the distribution function
$\det(I-K_{H,n})$ has an alternative representation as an $n\times n$
determinant.  Explicitly,
\[ F_n^{GUE}(s)=\det\left(\delta_{i,j}-\int_{s/\sqrt{2}}^\iy
\vp_i(x)\vp_j(x)\,dx\right)_{0\le i,j\le n-1}\]
where $\vp_j$ are the oscillator functions previously introduced.
This last representation implies that the $F_n^{GUE}$
are expressible in terms of elementary functions and
the error function with
increasing complexity for increasing values of $n$.
In the  simplest case, $n=1$,
 $F_1^{GUE}$ is the standard normal; a result easily anticipated from
the original formulation of the growth model.
The next simplest case is $n=2$,
\[ F_2^{GUE}(s)= {1\ov 4}-{1\ov 2\pi}\,
 e^{-s^2}-{1\ov 2^{3/2}\sqrt{\pi}}\,s \,e^{-s^2/2}+
{1\ov 2}(1-{1\ov\sqrt{2\pi}}s e^{-s^2/2})\,
\textrm{erf}(s/\sqrt{2})+{1\ov 4}\,
\,\textrm{erf}(s/\sqrt{2})^2.\]

The moments of $F_n^{GUE}$ are, of course,
\[ \mu_j(n):=\int_{-\iy}^\iy s^j f_n^{GUE}(s)\,ds, \> j=1,2,\ldots \]
where $f_n^{GUE}=dF_n^{GUE}/ds$.
For $1\le n \le 5$ we have,

\noindent First Moments:
\ba
\mu_1(1)&=&0,\\
\mu_1(2)&=&{2\ov\sqrt{\pi}} \approx 1.128379,\\
\mu_1(3)&=& {27\ov 8\sqrt{\pi}} \approx 1.904140,\\
\mu_1(4)&=&-{7\ov 48\sqrt{2}\,\pi^{3/2}}+{475\ov 128\sqrt{\pi}}+
{475\ov 64\,\pi^{3/2}}\,\arcsin(1/3) \approx 2.528113, \\
\mu_1(5)&=&{13715\ov 4096\sqrt{\pi}}-
{16975\ov 41472\sqrt{2}\,\pi^{3/2} }+{41145\ov
2048\,\pi^{3/2}}\,\arcsin(1/3) \approx 3.063268.
\ea
Second Moments:
\ba
\mu_2(1)&=&1,\\
\mu_2(2)&=&2,\\
\mu_2(3)&=&3+{9\sqrt{3}\ov 4\pi} \approx 4.240490,\\
\mu_2(4)&=&4+{16\ov \sqrt{3}\,\pi} \approx 6.940420,\\
\mu_2(5)&=&5-{155\sqrt{5}\ov 864\,\pi^2}+{2495\ov 108\sqrt{3}\,\pi}+
{499\ov 54\sqrt{3}\,\pi^2}\,\arcsin(1/4) \approx 1.977575.
\ea
Third Moments:
\ba
\mu_3(1)&=&0,\\
\mu_3(2)&=&{7\ov  \sqrt{\pi}}\approx 3.949327, \\
\mu_3(3)&=&{297\ov 16\sqrt{\pi}} \approx 10.472769, \\
\mu_3(4)&=& {333\ov 32\sqrt{2}\,\pi^{3/2}}+{7109\ov 256\,\pi^{1/2}}+
{7109\ov 128\,\pi^{3/2}}\,\arcsin(1/3) \approx 20.378309,\\
\mu_3(5)&=&{2595475\ov 82944\sqrt{2}\,\pi^{3/2}}+{259385\ov
8192\sqrt{\pi}}
+{778155\ov\ 4096 \,\pi^{3/2}} \approx 33.432221.
\ea
Fourth Moments:
\ba
\mu_4(1)&=& 3,\\
\mu_4(2)&=& 9, \\
\mu_4(3)&=& {19}+{33\sqrt{3}\ov 2\,\pi} \approx 28.096927, \\
\mu_4(4)&=& {33}+{496\ov 3\sqrt{3}\,\pi} \approx 63.384348, \\
\mu_4(5)&=& {51}+{7475\sqrt{5}\ov 1296\,\pi^2}+
{99575\ov 324\sqrt{3}\,\pi}+{99575\ov 162\sqrt{3}\,\pi^2}\,
\arcsin(1/4) \approx  117.872208.
\ea
\begin{table}
\begin{center}
{\small
\begin{tabular}{|r|cc|cc|c|c|}\hline
$n$ & $\mu$ & Approx & $\sigma^2$ & Approx
& $S$
 & $K$\\[0.5ex]
\hline
2 &  1.12838 & 1.251 & 0.72676  &  0.645 & 0.08465 & 0.01053\\
3 &  1.90414 & 1.989 & 0.61474 & 0.564 & 0.11862 & 0.02192 \\
4 &  2.52811 & 2.594 & 0.54907 & 0.512 & 0.13749 & 0.03042 \\
5 &  3.06327 & 3.118 & 0.50426 & 0.476 & 0.14972 & 0.03683 \\ \hline
6 &  3.53861 & 3.585 & 0.47101 & 0.448 & 0.15838 & 0.04184 \\
7 &  3.97026 & 4.011 & 0.44497 & 0.425 & 0.16490 & 0.04586 \\
8 &  4.36822 & 4.405 & 0.42379 & 0.407 & 0.17001 & 0.04917 \\
9 &  4.73920 & 4.772 & 0.40609 & 0.391 & 0.17414 & 0.05195 \\ \hline
\end{tabular}
}
\vspace{2ex}
\caption{\label{momentData}
The mean ($\mu$) and
the variance ($\sigma^2$) of $H_n^\iy$, $2\le n \le 9$,
($H_n^\iy$ has distribution function $F_n^{GUE}$) are compared with
the  approximations (\ref{airyApprox1})
and (\ref{airyApprox2}), respectively. Also displayed are the skewness
($S$) and excess kurtosis ($K$) of $H_n^{\iy}$.
  $H^\iy$ has  $S\approx 0.2241$ and
 $K\approx 0.0935$. }
\end{center}
\end{table}
Let $H_n^{\iy}$ denote the weak limit $m\ra\iy$, $n$ fixed, of
\[ {H-p\,m\ov \s \sqrt{m}}\, , \]
and $H^{\iy}$ the weak limit $m\ra\iy$, $n\ra\iy$,
$\al =n/m$  fixed, of
\[{1\ov v (3b)^{1/3}\, m^{1/3}} \left(H-c\,m\right).  \]
(Thus the distribution functions of
$H^\iy_n$ and $H^\iy$ are $F_n^{GUE}$ and
 $F_2$, respectively.)
For $\al\ra 0$, $c=p+2\s\sqrt{\al}+O(\al)$, $(3b)^{1/3}\sim
\s/\al^{1/6}$, and $v\sim 1$.
Proceeding heuristically,
\ba
H &\sim & c \, m + m^{1/3} v (3b)^{1/3}\, H^\iy \\
&\sim & p \, m + 2\s\sqrt{\al}\, m + m^{1/3} \s \al^{-1/6} \, H^\iy \\
&\sim & p \, m +\s \, m^{1/2} \left\{ 2\sqrt{n}+{H^\iy\ov
n^{1/6}}\right\}\, .
\ea
Thus we expect
\[ H_n^{\iy} \sim 2\sqrt{n}+{H^\iy\ov n^{1/6}}\, , \]
and hence
\bae E(H_n^\iy)&\approx &2\sqrt{n}+{E(H^\iy)\ov n^{1/6}},\>\>
E(H^\iy)=-1.77109\ldots,
\label{airyApprox1} \\
 \textrm{Var}(H_n^\iy)&\approx&
 {\textrm{Var}(H^\iy)\ov n^{1/3}},\>\>
\textrm{Var}(H^\iy)=0.8132\ldots.  \label{airyApprox2}\eae
 For $2\le n\le 9$,
 these approximations are compared with the exact moments in Table
\ref{momentData}. We also
 compute the skewness
and the excess kurtosis\footnote{The skewness of a random
variable $X$ is $E\left(({X-\mu\ov \s})^3\right)$ and the
excess kurtosis  is $E\left(({ X-\mu\ov \s})^4\right)-3$.
Here $\mu=E(X)$ and $\s^2=\textrm{Var}(X)$.} of $F_n^{GUE}$.

\begin{figure}
\hspace{2cm}\resizebox{10cm}{6cm}{\includegraphics{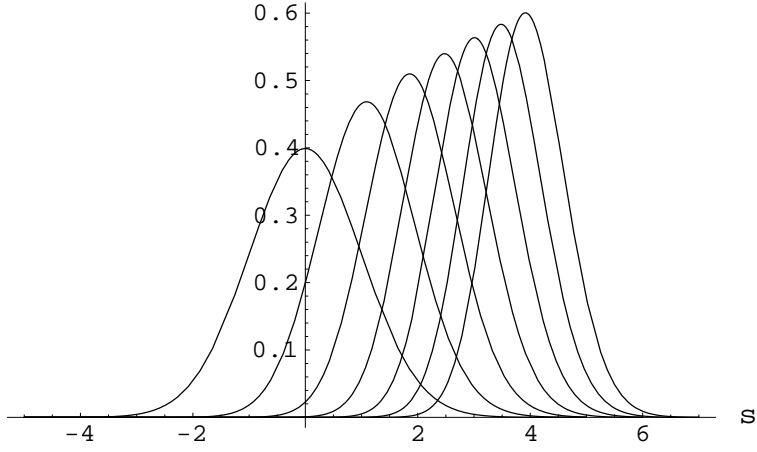}}
\caption{\label{densityFig}The densities $f_n^{GUE}$, $1\le n \le 7$.
For increasing values of $n$ the maximum of $f_n^{GUE}$ moves to
the right.}
\end{figure}

The densities $f_n^{GUE}$, $1\le n \le 7$, are graphed in
Fig.~\ref{densityFig}.

\section[Brownian Motion Representation in Finite $x$ Regime]{Brownian
Motion Representation  in the Finite $x$ GUE Regime}
Let $B(t)=(B_0(t), \dots, B_x(t))$
be the $(x+1)$-dimensional Brownian motion. Let $F$ be the
following functional on continuous functions $f=(f_0,\dots, f_x)$
from $[0,1]$ to $\textbf{R}^{x+1}$, which satisfy $f(0)=0$,
\[
F(f):=\max\left\{f_0(t_0)+f_1(t_1)-f_1(t_0)+\dots+f_x(t_x)-f_x(t_{x-1}):
0\le t_0\le t_1\le \dots \le t_x=1\right\}.
\]
Note that $F$ is continuous in  the $L^\infty$ metric.
Finally, let
\[ M_x:=F(B).\]

\textbf{Theorem.} For  $x\in\textbf{Z}_+$,
and $t\ra\iy$, we have
\[
\frac {h_t(x)-p t}{\sigma \sqrt t}\convd M_x,
\]
where $\sigma^2=p(1-p)$.

\textbf{Proof.} First, the path representation tells
us that we change $h_t(x)$ by at most a constant
if we only obey the  increasing property
within the same column. That is, we
change $L(x,t)$ to $L'(x,t)$, where
$L'(x,t)$ is the longest path $(x_i,t_i), i=1,\dots, k$,
of marked points such that
$0\le t_i-t_{i-1}-1$ if $x_i=x_{i-1}$, while
$0\le x_i-x_{i-1}\le t_i-t_{i-1}-1$ if $x_i\ne x_{i-1}$.
Thus, the first observation  is
$$
|L(x,t)-L'(x,t)|\le  x.
$$

Let $S_k^i$ equal the length of the longest increasing
sequence of points $(t,i)$, $0\le t \le k$.
Then
\be
L'(x,t)=\max\{S^0_{k_1}+S^1_{k_2}-S^1_{k_1}+\dots
+S^x_{k_x}-S^x_{k_{x-1}}:
0\le k_1\le k_2\le \dots\le k_x=t-x\}.
\label{LPrime} \ee
Now for every fixed $i$, $S_k^i$ is independent of $S_k^j$
for $j\ne i$. Let $X_k^i$ equal the indicator
of the event that $(i,k)$ is a marked point.
 Of course,
$S_k^i=\sum_{\ell=1}^k X_\ell^i$.

Let $S^i(\t)$, $\t\in\textbf{R}_+$, equal $S_k^i$ when $\t=k$ and be
obtained
by linear interpolation off the integers. Moreover, let
$\tS^i$ be the centered versions
$\tS^i(\t)=S^i(\t)-p \t$, and $\tS(\t)=(\tS_0(\t),\dots,
\tS_x(\t))$.
For $0\le\t\le1$, define
\[ X_t(\t):={\tS(t\t)\ov \s\sqrt{t}}\, ,\]
then the standard invariance
principle (see, e.g.\ \cite{durrett}, Ch.~7.) implies
that $X_t$ converges as $t\to\infty$
in distribution to the $(x+1)$-dimensional Brownian
motion $B$.

Now define
\be
L''(x,t)=\max\{\tS^0(t_0)+\tS^1(t_1)-\tS^1(t_0)+\dots +\tS^x(t_x)-
\tS^x(t_{x-1}) :
0\le t_0\le t_1\le \dots\le t_x=t\},
\label{LPrime2}
\ee
then
\[
\left|L'(x,t)-p t-L''(x,t)\right|\le 5x.
\]
(The linear interpolation
gives an error of at most four at each $t_i$ and we
incur an additional $x$
by replacing $t-x$ by $t$.) Note now that
(by making a substitution $t_i'=t_i/t$)
\[
\frac {L''(x,t)}{\sigma \sqrt t}=F\left(X_t\right).
\]
It follows (by continuity of $F$) that the theorem holds
with $L''(x,t)$ in place of $h_t(x)$, but this is clearly enough.

\textit{Remark 1}.
We should note that this theorem clearly holds in  more general
circumstances. For example, we could make every $\times$ count
an independent random number of jumps. We would get the same
theorem, with the only assumption that the said random
number has finite variance.

\textit{Remark 2}. The theorem also holds
for random words over an alphabet
with $n$ letters~\cite{tw5}, except that the $x+1$
Brownian motions are not independent, but they have to
sum  to 0, so the covariances
$\Gamma_{ij}$ equals $(n-1)/n^2$ when $i=j$ and
$-1/n^2$ otherwise. In the case of two equiprobable
letters, the limiting distribution of the centered and normalized
length of the longest weakly increasing subsequence
in a random word is equal to
the distribution of the random variable
\ba X &=&\max_{0\le t\le 1}\left(B_0(t)+B_1(1)-B_1(t)\right)\\
 &=& 2\max_{0\le t \le 1}\left(B_0(t)\right)-B_0(1)\\
&=&2M-N.
\ea
($M$ denotes the random variable
 $\max_{0\le t \le 1}B_0(t)$ and
$N$ denotes the random variable  $B_0(1)$.)
{}From the reflection principle it follows (see, e.g.\ pg.~395 in
\cite{durrett})
that the joint density of $(M,N)$ is
\[ f_{M,N}(m,n)=\sqrt{2\ov\pi}\,(2m-n)\, e^{-(2m-n)^2/2}, \>
\textrm{for}\>\>
m\ge 0,\,  m\ge n. \]
Thus the density of $X$ equals\footnote{C.~Grinstead, in
unpublished notes, also found a random walk interpretation of the
two-letter random word problem and used this to determine
the limiting distribution in this case.}
\[ f_X(x)=\int_0^x f_{M,N}(m,2m-x)\, dm =\sqrt{2\ov \pi}\, x^2
e^{-x^2/2}.
\]

\textit{Remark 3}. Limiting distribution of the centered and
normalized $h_t(1)$:
Here we have
\[
M_1=\max_{0\le t\le 1} (B_1(t)+ (B_2(1)-B_2(t)))=
\max_{0\le t\le 1} (B_1(t)- B_2(t)) + B_2(1)=M+N.
\]
(The random variables $M$ and $N$ are defined by the
last equality; and therefore, are not to be confused
with the random variables of the previous remark.)
Note that $N$ is standard normal. Since
 $(B_1-B_2)/\sqrt{2}$ is the standard Brownian
motion, $M$ equals, in
distribution,
$\sqrt{2}\,|N|$ again by
the reflection principle.   Even though $M$ and $N$ are
not independent,  $E(M_1)=\sqrt 2E(|N|)= 2/\sqrt\pi$. Moreover,
the conditional distribution of $N$ given
the entire path of $W:=B_1(t)- B_2(t)$, $0\le t\le 1$,
depends only on its final point $W(1)$.
Given this
final point $S$ equals $s$, the distribution is normal with mean $-s/2$
and variance $1/2$.
That is, if $\F_t$ is the Brownian filtration for $W$,
$S=W(1)$, then
\[ \pr\left(N\in dn\vert \F_1\right)=\pr(N\in dn\vert S=s)
={1\ov \pi}\,e^{-(x+s/2)^2}\,dn. \]
This makes it immediately possible
to compute the second moment of $M_1$, since
\[
E(MN)=E\left(E(MN|\F_1)\right)=E\left(M E(N|\F_1)\right)=-E(MS)/2=
-E\left((M/\sqrt 2)(S/\sqrt 2)\right)=-1/2,
\]
by a straightforward computation with the joint density above.
Therefore $E(M_1^2)=E(M^2)+E(N^2)+2E(MN)=2$.

In this way, the  density of $M_1$ is
\ba
f_{M_{1}}(x)&=&E\left(\pr\left(M_1=x|M,S\right)\right)\\
&=&\int f_{M,S}(m,s)\,\pr\left(M_1=x|M=m,S=s\right)\, dm ds\\
&=& \int_0^\infty  dm\int_{-\infty}^m
\frac 12 f_{M,N}(m/\sqrt 2, n/\sqrt 2)
\, \frac 1{\sqrt\pi} e^{-(x-m+n/2)^2}\, dn.
\ea
An explicit evaluation shows this last integral equals, as it must,
$f_2^{GUE}(x)$.

\textit{Remark 4}. Since the distribution function of $M_x$ equals
$F_{x+1}^{GUE}$, it follows from RMT~\cite{mehta} that we have the
alternative representation
\be
\pr\left(M_x\le s\right)=c_n \int_{-\iy}^{s}\cdots \int_{-\iy}^{s}
\Delta(x)^2\, e^{-{1\ov 2}\sum x_j^2}\,dx_1\cdots dx_n
\label{distrFnGUE}\ee
where
\[ \Delta(x)=\Delta(x_1,\ldots,x_n)=\prod_{1\le i<j\le n}(x_i-x_j)\]
is the Vandermonde determinant,
$c_n^{-1}=1!2!\cdots n! \, (2\pi)^{n/2}$,
and $n=x+1\ge 2$.  In the context of
Brownian motion, can one \textit{directly} prove
(\ref{distrFnGUE})?

\textit{Remark 5}. For connections between Brownian motion exit
times and random matrices, see Grabiner~\cite{grabiner}.

\textit{Remark 6}. The Brownian motion functional $M_x$ has appeared
previously in Glynn and Whitt \cite{glynn}, and consequently
(\ref{distrFnGUE}) provides an exact formula for the limiting
distribution of the departure time of the first $(x+1)$ customers from
$n$ single server queues.  Glynn and Whitt also consider the
case when $x=t^a$, $0<a<1$, and prove what would, in our setting,
be the following limit theorem
\[ \lim_{t\ra\iy} {h_t(x)-pt\ov \sigma\sqrt{tx}}=\alpha:=\lim_{x\ra\iy}
{M_x\ov \sqrt{x}},\]
with both limits in probability.
They conjectured that $\alpha=2$, and this was later proved by
Sepp\"al\"ainen~\cite{sepp2} via a
hydrodynamic limit for simple exclusion.  We note that our paper proves
that $\alpha=2$ as well, by a completely different route.  
 Namely, one only needs to apply the result (see, e.g.\ \cite{bai})
that the largest eigenvalue in the finite $n$ GUE scales as
$2\sqrt{n}$.\footnote{We thank Timo Sepp\"al\"ainen for bringing
this connection with queuing theory to our attention.}

\noindent\textbf{\large Acknowledgments}

This work was supported, in part,
by the National Science Foundation through grants
DMS--9703923, DMS--9802122 and
DMS--9732687.  In addition, the first author was supported
in part by the Republic of Slovenia's Ministry of Science,
grant number J1--8542--0101--97.  It is our pleasure to
acknowledge Iain Johnstone, Bruno Nachtergaele,
Timo Sepp\"al\"ainen
and Richard Stanley for helpful comments.  Finally, we wish
to thank both referees for their helpful comments and suggestions.

\end{document}